\documentclass[]{elsarticle}

\usepackage{graphicx, epsfig, amssymb, amsmath}
\usepackage{hyperref}
\usepackage[percent]{overpic}
\usepackage{color}
\usepackage{soul} 

\graphicspath{{./}{figs/}}

\newcommand{\REV}[1]{\textcolor{black}{#1}}
\newcommand{\R}{\mathbb R}

\newcommand{\Vd}{\mathbf V^{\textup{des}}}
\newcommand{\Vi}{\mathbf V^{\textup{int}}}
\newcommand{\Vo}{\mathbf V^{\textup{obs}}}
\newcommand{\V}{\mathbf V}
\newcommand{\vv}{\mathbf v}
\newcommand{\x}{\mathbf x}
\newcommand{\y}{\mathbf y}
\newcommand{\Sr}{\mathcal{S}_r}
\newcommand{\Cr}{C^{\textup{rep}}}
\newcommand{\E}{\mathcal{E}}
\renewcommand{\O}{\mathcal{O}}

\newcommand{\rhoin}{\rho_{\textup{in}}}

\newtheorem{rmk}{Remark}[section]

\journal{Applied Mathematical Modelling}

\begin{document}

\begin{frontmatter}

\title{Handling Obstacles in Pedestrian Simulations: \\ Models and Optimization}

\author[IAC]{Emiliano Cristiani}
\address[IAC]{Istituto per le Applicazioni del Calcolo, Consiglio Nazionale delle Ricerche, Rome, Italy.
\texttt{\{e.cristiani,d.peri\}@iac.cnr.it}}

\author[IAC]{Daniele Peri}


\begin{abstract}
In this paper we are concerned with the simulation of crowds in built environments, where obstacles play a role in the dynamics and in the interactions among pedestrians. First of all, we review the state-of-the-art of the techniques for handling obstacles in numerical simulations. Then, we introduce a new modelling technique which guarantees both impermeability and opacity of the obstacles, and does not require \textit{ad hoc} runtime interventions to avoid collisions. Most important, we solve a complex optimization problem by means of the Particle Swarm Optimization method in order to exploit the so-called Braess's paradox. More precisely, we reduce the evacuation time from a room by adding in the walking area multiple obstacles optimally placed and shaped.
\end{abstract}


\end{frontmatter}


\section{Introduction}\label{sec:intro}
In this paper we are concerned with the simulation of crowds in built environments, where obstacles play a role in the dynamics and in the interactions among pedestrians. Although obstacles are commonly included in numerical simulations, in many cases no special attention is given to the (numerous) related issues. Usually the method used to handle obstacles is poorly or not at all described, sometimes  obstacles are processed as normal boundary conditions or even as frozen pedestrians. We think instead that this matter deserves more attention, especially considering its importance and impact in the simulations. In particular, we refer to simulations which investigate the well known Braess's paradox \cite{braess2005TS, hughes2003ARFM}, which states that an additional obstacle or constraint can improve global dynamics.
Note that in the case of crowds, placing an additional obstacle may be intuitively seen as a bad idea. Nevertheless, a well placed obstacle can decrease the internal pressure among pedestrians and break symmetries in front of an exit, resulting in a faster outflow, see, e.g., \cite[\S 6.3]{albi1504.04064}. 
The obstacles can also facilitate the evacuation ensuring that all the available exits are equally used, see, e.g.,  \cite{cristiani2015SIAP}. 
It is then crucial that interactions between crowd and obstacles are correctly handled.

\medskip

As an introduction to the field of pedestrian modelling, we refer the reader to the \REV{survey papers \cite{bellomo2011SR, duives2013TRC, ronchi2016} and the books \cite{cristiani2014book, kachroo2008book}. }
Here we consider a macroscopic description of the crowd based on a two-dimensional first-order nonlocal conservation law already considered in \cite{cristiani2010pareschibook, cristiani2011MMS, cristiani2014book, cristiani2015SIAP, piccoli2009CMT, piccoli2011ARMA}. The peculiarity of this formulation is that no fundamental diagram is involved. Moreover, it is proven to be able to catch several self-organizing phenomena actually observed in crowds (see \cite[\S 1.1.2]{cristiani2014book} for some examples). However, the model is not naturally endowed with the capability of handling obstacles, which must be added as an independent feature.

\medskip 

To begin with, we review all the existing (in the knowledge of the authors) techniques employed to handle obstacles in pedestrian simulations, pointing to the main related references. Advantages and drawbacks of each method are discussed. This section can be useful to researchers who are entering the field.

Second, we introduce a new modelling technique to handle obstacles which is both simple to implement and effective for realistic simulations, thus giving a valid alternative to the existing techniques. The new method enforces that \textit{obstacles are impermeable} (regardless of the model parameters) and \textit{collisions are avoided} without \textit{ad hoc} runtime interventions. Pedestrians bypass the obstacle smoothly, \REV{are not trapped} in bays, and, especially, they do not see through the obstacles, which are assumed to be \textit{opaque}. This avoids unrealistic effects around thin walls.

Third, we tackle a challenging problem which consists in finding the optimal position and shape of some obstacles so that the global dynamics of pedestrians are improved (e.g., by a reduction of the evacuation time from a room). 
Several papers investigate numerically the effectiveness of the Braess's paradox by means of both microscopic models (e.g., Helbing's social force model) and macroscopic models, reporting the effect of additional obstacles manually placed in the walking area. See, among others, \cite{escobar2003LNCS, frank2011PA, helbing2005TS, hughes2002TRB, matsuoka2015, twarogowska2014AMM}. 
In this paper, instead, we follow the lines of \cite{cristiani2015SIAP, johansson2007, shukla2009}, where an optimization algorithm is used. The main novelty here is that we consider \textit{multiple free-shaped} obstacles, optimally placed and shaped by means of a modified Particle Swarm Optimization (PSO) method. The PSO is expected to give better results than genetic algorithms \cite{johansson2007, shukla2009} and random compass search \cite{cristiani2015SIAP}.

\REV{
Finally, let us stress the importance of pedestrians guidance in the context of emergency situations. Such a circumstances -- although rarely characterized by actual state of \textit{panic} (intended as irrational, competitive, and non altruistic behaviour) \cite{cocking2014JCASP, fahy2012FM} -- involve fast decision process which can lead to non optimal route choice. For example, it is known that heading toward the nearest exit is not always the best evacuation strategy \cite{lin2015NHM}. Moreover, people could be influenced by others (social influence or herding effect) and tend to follow people who show definiteness \cite{albi1504.04064, kinateder2014TRF}. This urges us to design environments where \textit{natural} (intended as actually observed) behaviour ideally coincides with the \textit{optimal} (with respect to some criterion to be defined) behaviour. This can be achieved by means of suitably placed obstacles, as already investigated in \cite{cristiani2015SIAP}.
}

\medskip 

\noindent \textit{Paper organization.} 
In section \ref{sec:review} we review the existing techniques used to handle obstacles in numerical simulations.
In section \ref{sec:model} we present the new method to include obstacles and a first numerical test to show its main features.
In section \ref{sec:shapeoptimization} an outline of the optimization algorithm, together with the description of the techniques for obstacle parameterization and management, is presented. 
In section \ref{sec:preliminarynumericalinvestigation} some preliminary tests put some light on the convergence properties and the numerical accuracy of the model.
Section \ref{sec:numericaltests} contains the results of the optimization of a \REV{simple} environment.

\section{A brief review of obstacles' handling techniques}\label{sec:review}
In this section we briefly review the most common techniques used to deal with obstacles in the literature about pedestrian modelling. A general-purpose survey of pedestrian models can be found in \cite{cristiani2014book}. Here we just recall the two main ingredients of pedestrian models: (i) a \textit{desired velocity} which steer pedestrians toward a (common) target, for example an exit door. This is the velocity field people would follow if they were alone in the walking area; (ii) a \textit{Repulsion (social) force} exerted by pedestrians themselves, which accounts for the tendency of people to stay away from crowded regions and avoid collisions.

\medskip

\noindent $\bullet$ {\bf Repulsive obstacles.} 
One of the most common methods used to manage obstacles is obtained assuming that they generate a repulsive (social) force, exactly as pedestrians themselves do. In other words, obstacles are teated as frozen pedestrians. In this way one can use a repulsion function of the same kind to model both the interactions with group mates and with obstacles. \REV{This method is extensively used in microscopic models, see, e.g., \cite{colombi2016CAM, frank2011PA, helbing1995PRE, loehner2010AMM, okazaki1979TAIJa, okazaki1979TAIJb, okazaki1979TAIJc} and also in macroscopic and multiscale models, see, e.g., \cite{colombo2012M3AS, coscia2008M3AS, etikyala2014M3AS, piccoli2009CMT}, with or without the pre-evaluation of the distance-to-obstacle function.} The main drawback of this approach is that it is quite difficult to tune the strength of the repulsion force in such a way that the resulting behaviour is both admissible and realistic. Indeed, if the force is too small there is the risk that pedestrians enter the obstacles, while if it is too large pedestrians bypass the obstacles excessively far away. The paper \cite{colombo2012M3AS} proposes a method to tune automatically the strength of the repulsion. From the computational point of view, it is useful to note that interactions with obstacles must be computed continuously during the simulation.

\medskip

\noindent $\bullet$ {\bf Cut off of the velocity field.}
Another easy method to deal with obstacles is obtained \REV{by} computing the velocity field first neglecting the presence of the obstacles, then nullifying the component of the velocity vector which points inside the obstacle. This method is used in, e.g., \cite{albi1504.04064, cristiani2011MMS, cristiani2015SIAP}. Again, handling obstacles in this way is expensive from the computational point of view, since interactions with obstacles must be checked continuously during the computation. Moreover, one must be sure that pedestrians do not stop walking completely because both components of the velocity vector vanish. This can happen around corners, stair-shaped obstacles and when obstacles are very close to each other (i.e.\ the distance is comparable with the spatial resolution of the numerical grid).

\medskip

\noindent $\bullet$ {\bf Rational turnaround.}
In more sophisticated models which take into account the rationality and predictive ability of pedestrians, obstacles can be managed including them into the decisional process. For example, in the Hughes's model \cite{hughes2002TRB} the pedestrians move, at each time, along the fastest path toward the target, considering that crowded regions slow down the walking speed. In this framework, obstacles are easily included assuming that inside them the speed is null, so that the computation of the fastest path will circumvent them automatically. 
In \cite{carrillo2016M3AS} it is assumed that obstacles have a ``zone of influence''. This is translated by the fact that the admissible maximum speed is reduced as one approaches the obstacle, and vanishes inside it. 
In \cite{foderaro2014A} an hybrid approach is proposed: obstacles do generate a repulsive force, then this force contributes to define a cost functional which is minimized in a separate procedure. In the context of Cellular Automata models, rationality is often obtained by means of a \emph{floor field}, which affects the matrix of transition probabilities. 
In \cite{varas2007PA} the floor field is used to define a desired velocity to a target in such a way that pedestrians reach the target moving from a cell with large field to a cell with smaller field. Assigning a very large value to the field inside the obstacles makes pedestrians avoid naturally prohibited cells.

\medskip

\noindent $\bullet$ {\bf Collision course manoeuvres.}
Some crowd models assume that pedestrians can forecast the trajectory of the others. Then, each pedestrian checks continuously if he/she is on a collision course with some others. If this is the case, he/she responds in order to avert the danger of collision, see e.g., \cite{karamouzas2014PRL}. This can be done both with group mates and obstacles. For example, in \cite{agnelli2015M3AS} if a pedestrian is pointing toward an obstacle, his/her trajectory is modified steering toward the tangent to the obstacle.

\medskip

\noindent $\bullet$ {\bf No trespassing.}
In Cellular Automata models without floor field, obstacles are managed easily avoiding pedestrians to occupy cells representing obstacles. The prohibition is obtained nullifying the probability of the pedestrian to move in that cell, see, e.g., \cite{kirchner2002PA}.


\section{A new model with obstacle management}\label{sec:model}
Let us assume that pedestrians are free to move in a walking area $\Omega\subset\R^2$. 
People either enter the domain $\Omega$ through a door or they are assumed to be already in at the initial time. We denote by $\E$ the exit doors, which coincide with the target of pedestrians. Obstacles are denoted by $\O$.

\subsection{The basic model}\label{sec:basicmodel}
Following the lines of \cite{cristiani2015SIAP}, we describe the evolution of the average density $\rho:\Omega\times\R^+\to\R^+$ of the crowd by means of the following first-order nonlinear nonlocal conservation law
\begin{equation}\label{CL}
\left\{
\begin{array}{ll}
\rho_t(\x,t)+\nabla\cdot\Big(\rho(\x,t)\V[\rho(\cdot,t)](\x)\Big)=0, & \x\in\Omega\backslash(\O\cup\E), \quad t>0,\\ [2mm]
\rho(\x,0)=\bar\rho(\x),  & \x\in\Omega\backslash(\O\cup\E), \\ [2mm]
\rho(\x,t)=0,  & \x\in(\O\cup\E),\quad t>0,
\end{array}
\right.
\end{equation}
where the velocity field $\V:\Omega\to\Omega$ is split in two parts:
\begin{equation}\label{V_no_obs}
\V[\rho(\cdot,t)](\x):=\Vd(\x)+\Vi[\rho(\cdot,t)](\x).
\end{equation}
The \textit{desired velocity} $\Vd$ corresponds to the velocity each pedestrian would have if he/she was alone in the domain. It points toward the nearest exit along the shortest path (it does not always coincide with a straight line because of the obstacles) and does not change in time. Resorting to the Bellman's dynamic programming principle and the optimal control theory based on Hamilton-Jacobi equations \cite{bardibook}, $\Vd$ is computed as 
\begin{equation}\label{Vd}
\Vd(\x):=-\frac{\nabla\varphi(\x)}{|\nabla\varphi(\x)|}
\end{equation}
where $\varphi$ is the unique viscosity solution to the eikonal equation
\begin{equation}\label{eikonal}
\left\{
\begin{array}{ll}
|\nabla\varphi(\x)|=1, & \x\in\Omega\backslash(\E\cup\O), \\ [2mm]
\varphi(\x)=0,  & \x\in\E, \\ [2mm]
\varphi(\x)=+\infty,  & \x\in\O.
\end{array}
\right.
\end{equation}
This PDE is solved once, before the simulation starts. It is useful to note that $\varphi$ is the distance function from the exits $\E$ and that $\Vd$ never points inside obstacles by construction. Moreover, $\Vd(\x)$ is a unit vector.
 
The \textit{interaction velocity} instead is given by the following nonlocal term
\begin{equation}\label{Vint}
\Vi[\rho(\cdot,t)](\x):=\iint_{\Sr(x)}\frac{\y-\x}{|\y-\x|}\frac{-\Cr}{|\y-\x|}\rho(\y,t)d \y,\quad \x\in\Omega,
\end{equation}
where $\Cr>0$ is a positive constant (repulsion strength), and the sensory region $\Sr$ is defined as
\begin{equation}\label{S}
\Sr(\x):=\left\{\y\in\Omega\ :\ \varepsilon<|\y-\x|<r \text{ and } \frac{\y-\x}{|\y-\x|}\cdot \Vd>0\right\},\quad r>\varepsilon>0.
\end{equation}
The interaction velocity (\ref{Vint}) makes people move away from crowded region and its modulus is inversely proportional to the distance from the others. The sensory region \eqref{S} instead is a circular region in front of the pedestrian, assuming that his/her head is aligned with the desired velocity. The size of the sensory region is ruled by the parameter $r$, while the (small) parameter $\varepsilon$ is used to avoid singularities. Note that no fundamental diagram appears here and the density $\rho$ is no \textit{a priori} bounded by a constant (even if it cannot converge to a Dirac delta if $\varepsilon$ is sufficiently small).

Boundary and inflow conditions for (\ref{CL}) will be discussed afterwards together with the obstacles.
\subsection{Avoiding obstacles}
The dynamics described in section \ref{sec:basicmodel} need a correction in order to deal with obstacles. Indeed, while $\Vd$ does not point inside obstacles by construction, the interaction velocity $\Vi$ could do so (pedestrians could be ``pushed'' inside the obstacles by group mates). 
To fix the problem, let us first compute offline the distance function $\Phi=\Phi(x)$ from the obstacles $\O$, solving another eikonal equation
\begin{equation}\label{eikonal_obs}
\left\{
\begin{array}{ll}
|\nabla\Phi(\x)|=1, & \x\in(\Omega\backslash\O)\cup\E, \\ [2mm]
\Phi(\x)=0,  & \x\in\O.
\end{array}
\right.
\end{equation}

Now, let $\lambda=\lambda(\Phi)$ be any decreasing function of the distance from the obstacles such that $\lambda(0)=1$ and $\lambda(\Phi)=0$, for $\Phi>M>0$. We define the velocity field $\Vo$ in $\Omega$ as
\begin{equation}\label{Vobs}
\Vo[\rho](\x):=
\underbrace{\Big(1+\Vd(\x)\cdot\Vi[\rho](\x)\Big)^+}_{:=s[\rho](\x)}\ \Vd(\x),
\end{equation}
and then we finally correct the total velocity $\V$ in \eqref{CL}-\eqref{V_no_obs} with the following one
\begin{equation}\label{V}
\hat \V[\rho](\x):=\lambda(\Phi(\x))\Vo[\rho](\x)+(1-\lambda(\Phi(\x)))(\Vd(\x)+\Vi[\rho](\x)).
\end{equation}
The special construction of $\Vo$ guarantees that it never points inside obstacles since it is parallel to $\Vd$, which does not. The modulus of $\Vo$ is given by $s$, which acts as the classical relationship between speed and density (cf.\ the fundamental diagram). It is positive or zero, meaning that pedestrians can only decrease their distance from the target or leave it unchanged. It is then clear that the final velocity $\hat\V$ coincide with the original $\V$ sufficiently far from the obstacles, and it is gradually substituted by $\Vo$ as one approaches the obstacles.

\begin{rmk}\label{rem:advantages}
Our approach is quite advantageous from the computational point of view. First of all, one can compute and store in advance $\varphi$, $\Vd$, $\Phi$, and the points belonging to $\Sr$, because they do not depend on $\rho$. Second, pedestrians never enter the obstacles by construction, so there is no need to check every time if a correction of the trajectory is required. This features allow to save a huge amount of CPU time and pave the way to optimization algorithm, which usually need to run the simulation several times. This is what we do in section \ref{sec:PSO}.
\end{rmk}
\subsection{Handling opaque obstacles}\label{sec:handlingopaqueobstacles}
It is quite natural assuming that pedestrians cannot see through the walls and, in particular, that they are not at all affected by group mates beyond the obstacles. Nevertheless, only a few papers deal with the issue of opaque obstacles, see, e.g., \cite[Appendix A]{colombo2012M3AS}. Usually pedestrians are assumed to interact with any other group mate in the domain, even if it is hidden by an obstacle.

\begin{figure}[h!]
\begin{center}
\begin{overpic}
[width=0.28\textwidth]{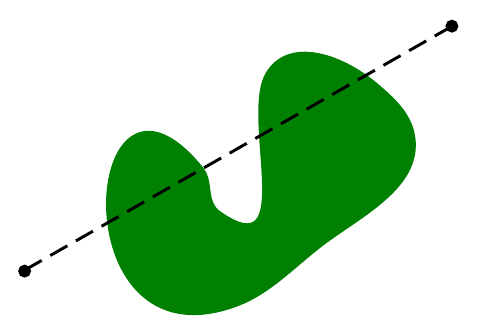}
\put(98,68){$\y$}
\put(2,2){$\x$}
\put(60,20){$\O$}
\end{overpic} 
\end{center}
\caption{Points $\x$ and $\y$ are not visible from each other.}
\label{fig:opaqueobstacles}
\end{figure}

In this paper we consider opaque obstacles, implementing vision obstruction also at the numerical level. A pedestrian located at the point $\x\in\R^2$ interacts with a pedestrian located at point $\y\in\R^2$ if the line joining $\x$ with $\y$ does not intersect any obstacle, see figure \ref{fig:opaqueobstacles}. 
At the discrete level, this can been easily implemented: once the domain $\Omega$ is discretized by a mesh, we divide the line joining $\x$ and $\y$ in $N_s$ evenly spaced points in such a way that two consecutive points lie in the same cell or in two adjacent cells.
Then one can check if any point falls in a cell labelled as obstacle and, if in this case, the point $\y$ is defined as ``non visible'' from $\x$ and consequently removed from $\Sr$. Again, this operation can be done offline before simulation starts, see Remark \ref{rem:advantages}.

\REV{
\begin{rmk}\label{rem:limitations}
Despite the fact that obstacles are opaque and people are not affected by group mates beyond the obstacles, we still assume that people are able to find the shortest path to the target (among all the possible routes, even if they are not visible). 
If the environment is instead assumed to be unfamiliar, this will not be possible. In this case people will take dynamic strategic decisions according to their environmental awareness, which increases due to new information acquired during the journey, see, e.g., \cite{colombi2016CAM}. 
Even if people are familiar with the environment, they might choose other paths than the shortest one: in Hughes-like models, the fastest path is continuously recomputed taking into account the current distribution of pedestrians (assumed to be known everywhere \cite{hughes2002TRB} or not \cite{carrillo2016M3AS}).  
Following the mean-field game approach (also known as predictive dynamic user equilibrium), people react both to current and future pedestrian distribution, being able to anticipate exactly future dynamics of the crowd \cite{cristiani2015SIAP, hoogendoorn2004TRB, lachapelle2011TRB}. 
Finally, let us mention that people may also be influenced by the others' choice (the so-called social influence or herding effect), see, e.g., \cite{albi1504.04064, kinateder2014TRF}.
\end{rmk}
}
\subsection{Numerical approximation}\label{sec:numerics}
The conservation law \eqref{CL} is discretized by means of the scheme firstly proposed in \cite{piccoli2011ARMA} and then used extensively in, e.g., \cite{cristiani2011MMS, cristiani2014book}. It is a two-dimensional first-order reasonably fast conservative scheme, which has been proved to converge to a weak solution of \eqref{CL}, see \cite{piccoli2013AAM, TosinFrasca}. It describes adequately the main features of pedestrian flow, including merging and splitting, although it exhibits a non-negligible numerical diffusion (its one-dimensional version coincides with the classical upwind scheme). 

The eikonal equations \eqref{eikonal} and \eqref{eikonal_obs} are discretized by means of an iterative first-order semi-Lagrangian scheme. The interested reader can find a complete introduction to the topic in the recent book \cite{falconebook} (see also \cite{bardibook}). The present scheme is described (with complete references) in \cite{cacace2014SISC}. The reconstruction of the values of the solution at non-mesh points is obtained by means of a three-point linear interpolation. The Fast Sweeping technique (see again \cite{cacace2014SISC} for explanations and references) is used to speed up the convergence.

The spatial domain is $\Omega=[0,1]^2$, discretized by means of a Cartesian grid. If not otherwise stated, we employ $100\times 100$ grid nodes. The time step is equal to the 90\% of the maximum value allowed by the CFL condition. The function $\lambda$ decrease linearly to 0 within a distance of 0.03 from the boundary of the obstacles.
\subsection{A first numerical example}
\begin{figure}[t!]
\begin{center}
\begin{tabular}{lr}
\begin{overpic}
[width=0.43\textwidth]{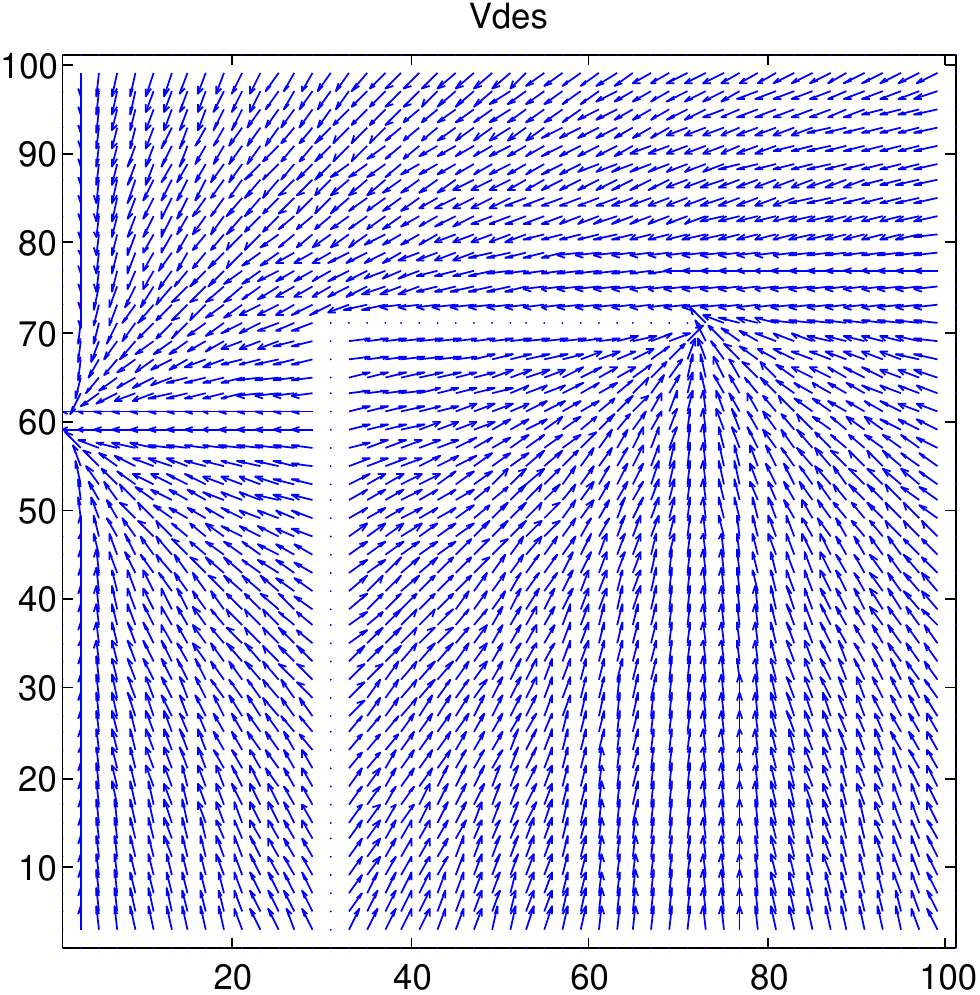}
\put(-4,55){$\E$}\put(36,30){$\O$}
\end{overpic}
&
\begin{overpic}
[width=0.48\textwidth]{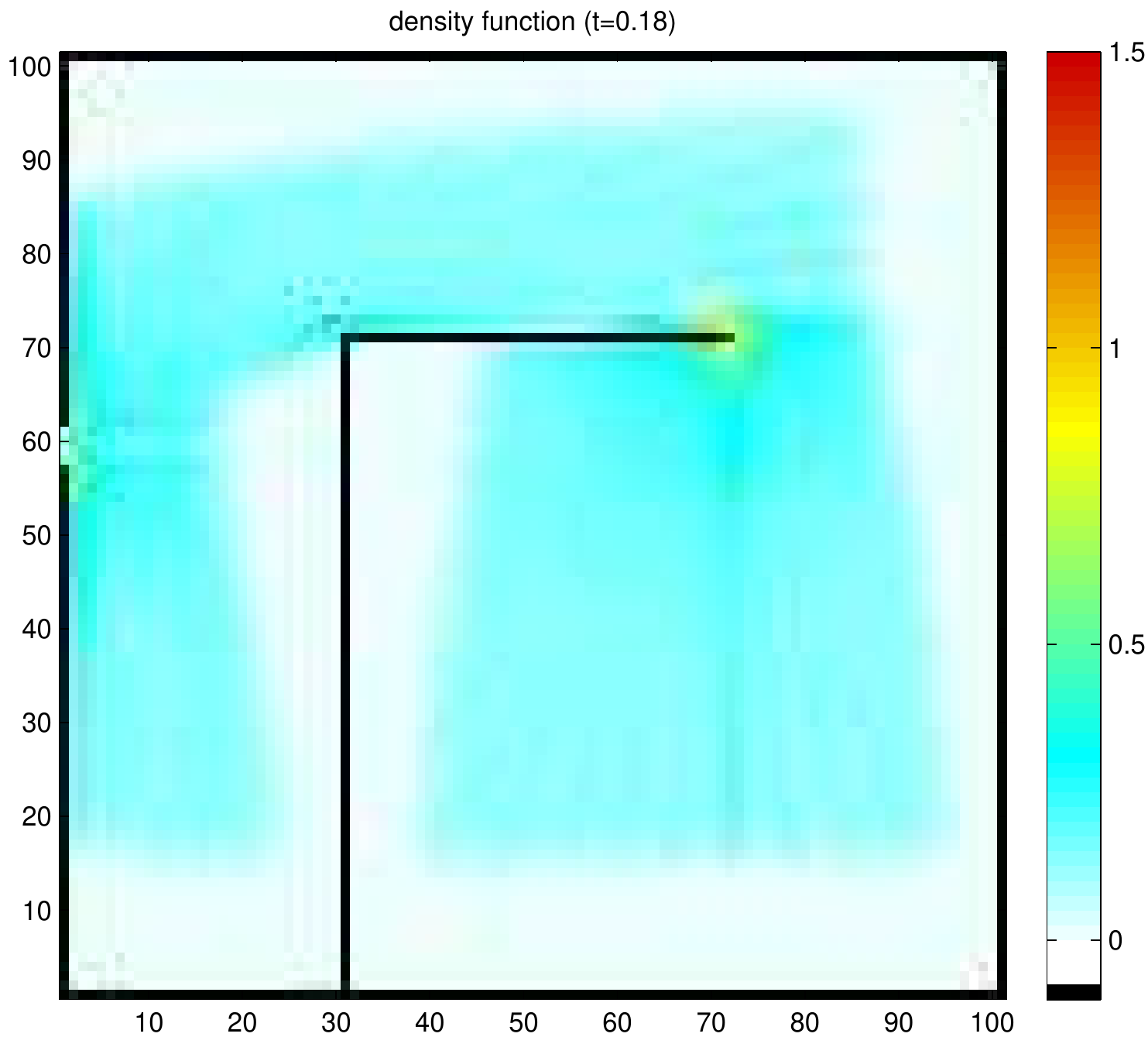}
\end{overpic}
\\
\begin{overpic}
[width=0.48\textwidth]{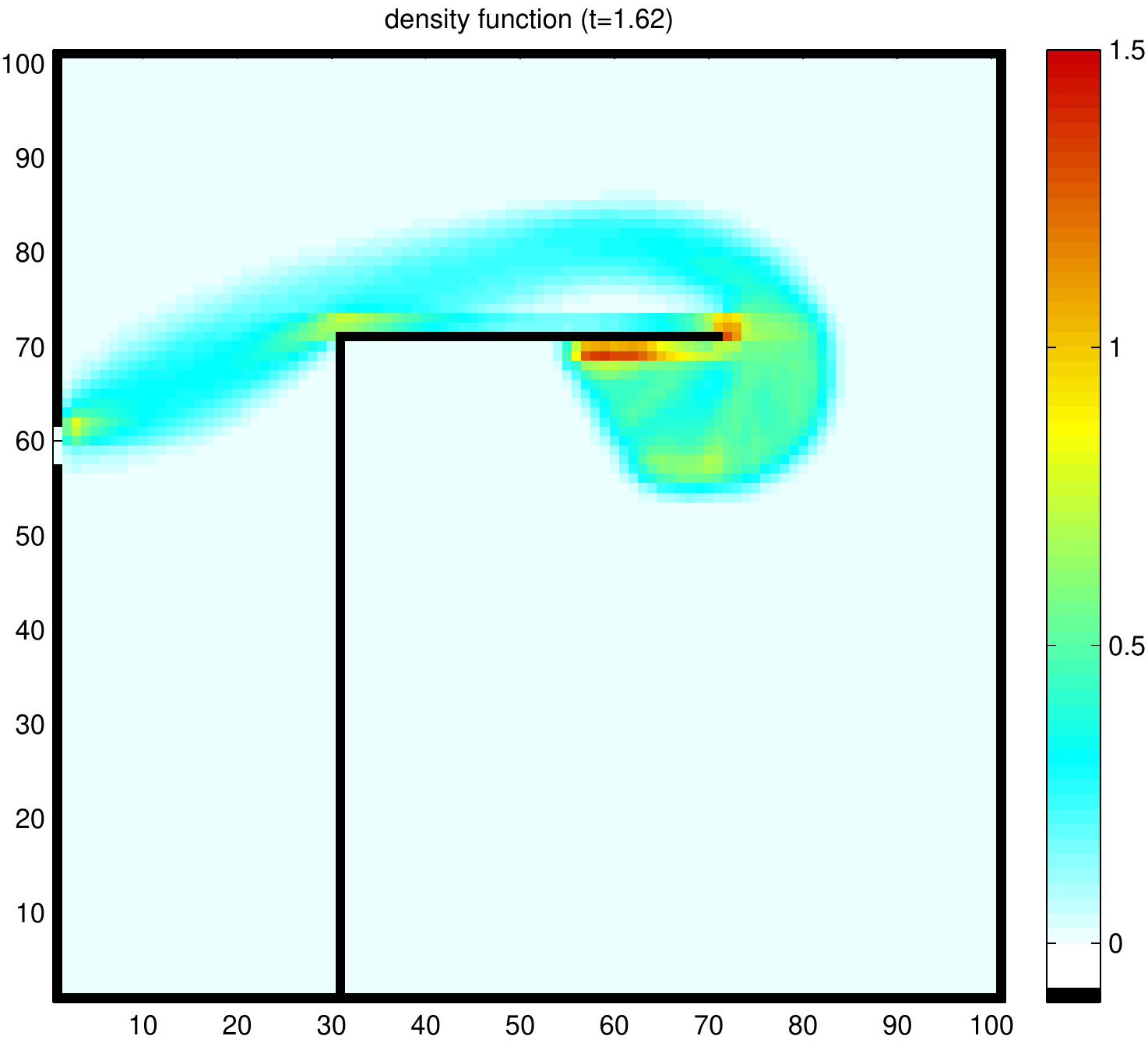}
\put(-4,50){$\E$}\put(32,30){$\O$}
\end{overpic}
&
\begin{overpic}
[width=0.48\textwidth]{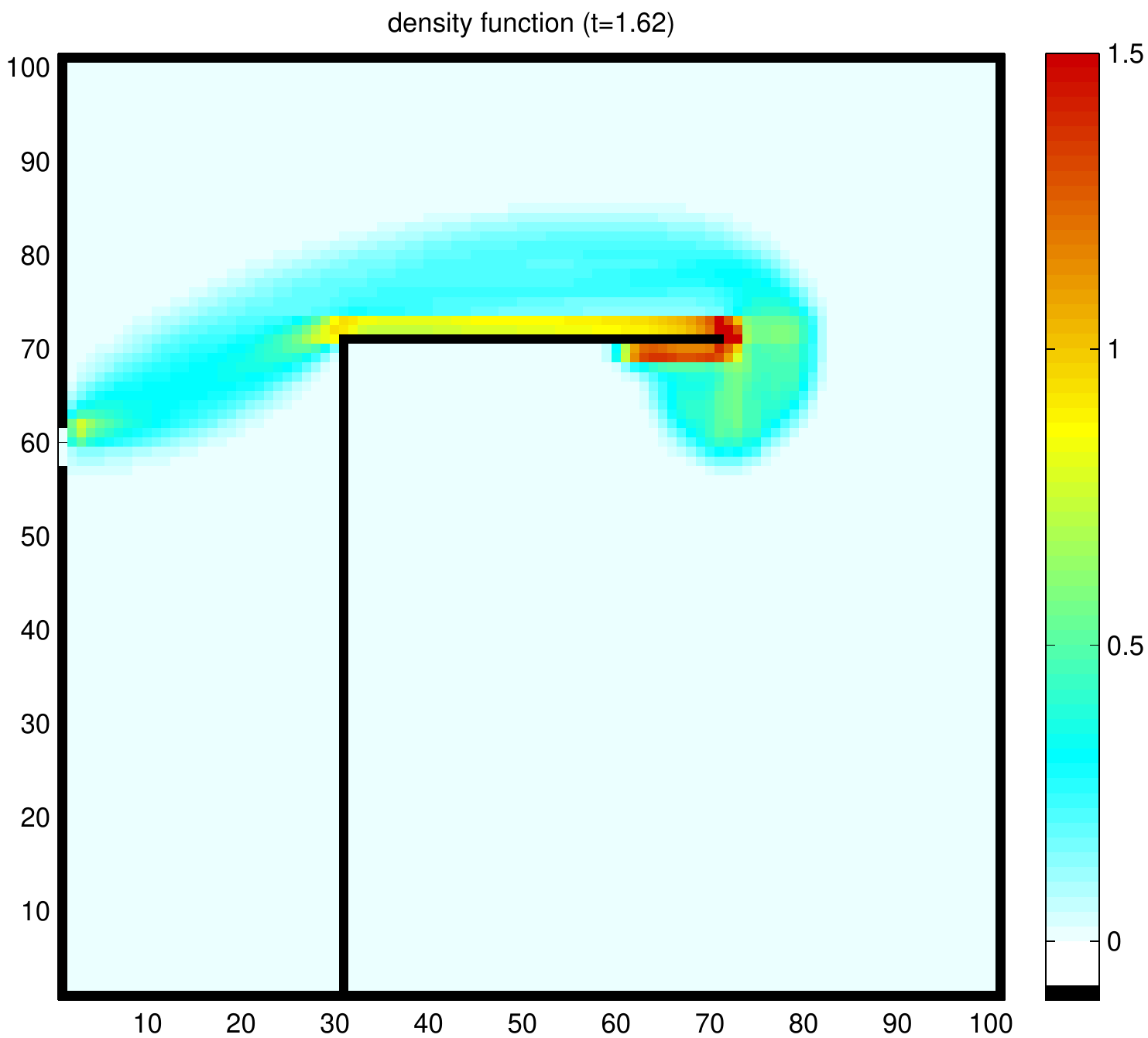}
\end{overpic}
\end{tabular}
\end{center}
\caption{Desired velocity field $\Vd$ (top-left), pedestrian density after a small time (top-right), pedestrian density at $t=1.62$ with transparent obstacle (bottom-left), pedestrian density at $t=1.62$ with opaque obstacle (bottom-right).}
\label{fig:vedononvedo}
\end{figure}

We present here a first numerical test in order to enlighten the main features of the proposed model, in particular the difference between transparent and opaque obstacles. 
We place one exit $\E$ (target) on the left side and a L-shaped obstacle $\O$ in the interior, see figure \ref{fig:vedononvedo}. At initial time, people are uniformly distributed in $\Omega\backslash(\E\cup\O)$. The parameters of the simulation are: $\Cr=6$, $r=0.16$, $\bar\rho\equiv 0.1$.  
In the first case, people interact with any group mate in $\Sr$. In the second case, obstacles are opaque and people do not interact with group mates located on the other side of the obstacle, even if they are in $\Sr$. 
Figure \ref{fig:vedononvedo} shows the desired velocity field $\Vd$ and the results of the simulations. 
The difference between the two cases is quite visible on the top of the horizontal part of the obstacle. In the first case, people walking leftward feel the others below, who still have to turn around the corner, and move up, being repulsed from them. This clearly does not happen in the second case. The final evacuation time is 4.18 in the first case and 3.80 in the second case. 


\section{Shape optimization framework}\label{sec:shapeoptimization}
Two are the main elements of the optimization framework: a parametric model for the
obstacle description, able to express the size, shape and position of an obstacle by a
limited number of parameters, and an optimization algorithm, able to find rapidly the
best design vector, that is, the more convenient obstacle configuration with respect to a given cost criterion. In the following, we will focus on \textit{evacuation problems}, meaning that the cost to be minimized will be the evacuation time, defined as the first time $t^*$ such that the total mass $\iint_\Omega\rho(\x,t^*)d\x<10^{-5}$.

\subsection{Parametric model}

There are a number of main features that are necessary for a parametric model to make
it suitable for the application to design optimization problem. Firstly, the number of
design parameters should be small: this is required in order to reduce the complexity
of the optimization problem, that increases with the dimension of the solution space. With
this small number of parameters, a wide variety of shapes should be described: the range
of shapes should be as large as possible, without an intrinsic preference for a certain
class of shapes. A bijection between the design vector and the shape of the obstacle
is also recommended, in order to prevent a cyclical behaviour of the objective function.
For this specific problem, convexity of the obstacle is also important, in order not to
produce some niches where the pedestrian may be trapped into. This feature is unnecessary
for the adopted mathematical model, since it is able to prevent this situation. However, 
it is preferable to preserve convexity in order to be more realistic.

In this application, an obstacle is described by using a closed \textit{B\'ezier curve}. In a parametric
form, the B\'ezier curve is expressed as
\[
 {\bf B}(t) = \sum_{i=0}^n \binom{n}{i}{\bf P}_{i}(1-t)^{n-i}t^i, \qquad t \in [0,1],
\]
where ${\bf P}_0,\ldots,{\bf P}_n$ are $n+1$ given points in $\R^2$.
The outcome is a curve ${\bf B} \in \R^2$.

The curve is obtained in parametric form, using one equation for each coordinate. 
The control points $\{{\bf P}_i\}_i$ form the so called \textit{B\'ezier polygon}. 
The passage through the extreme points ${\bf P}_0$, ${\bf P}_n$ is guaranteed for the B\'ezier curve. 
Then we easily obtain a closed curve by choosing ${\bf P}_n={\bf P}_0$. 
Moreover, we guarantee that the curve is smooth at the closure point by choosing ${\bf P}_{n-1}$ in such a way that tangent is continuous (avoiding a cusp). 
This point must be in line with ${\bf P}_0$ and ${\bf P}_1$,
in the opposite direction with respect to ${\bf P}_1$. 
This is possible because the tangent at the extremal points
of the B\'ezier curve coincides with the direction of the two adjacent control points. 
Summarizing, in order to produce a smooth closed curve, among the $n+1$ control points of the B\'ezier curve, only ${\bf P}_0,\ldots,{\bf P}_{n-2}$ need to be provided, resulting in $2(n-1)$ degrees of freedom.

The curve is contained inside the convex hull of the B\'ezier polygon, since it is obtained as an envelope of line segments connecting the edges of the polygon. 
The degree of the B\'ezier curve is connected with the number of control points, and a small number
of control points is avoiding too complex shapes.

A further step is needed in order to discharge self-intersecting curves and concave shapes in general.
Since, in practice, the mathematical model is requiring to check each grid point to
be an inner point or not (see section \ref{sec:handlingopaqueobstacles}), 
this check is performed by tracing the sign of the cross product
between the vector pointing from a grid point to a generic point on the B\'ezier curve
and the tangent vector at the same point of the B\'ezier curve. If the sign of the cross
product is constant along the entire B\'ezier curve, the grid point is considered an inner
point, otherwise it is considered external. The obstacle is finally implemented as the set of inner
points of the computational grid detected by using this technique. This approach has two advantages: 
if the B\'ezier curve is self-intersecting, no inner point is detected. Secondly, the inner points are
reduced to a convex set, and the points generating the concave part of the shape are automatically
discharged.

\begin{figure}[h]
\begin{center}
\includegraphics[width=0.32\textwidth]{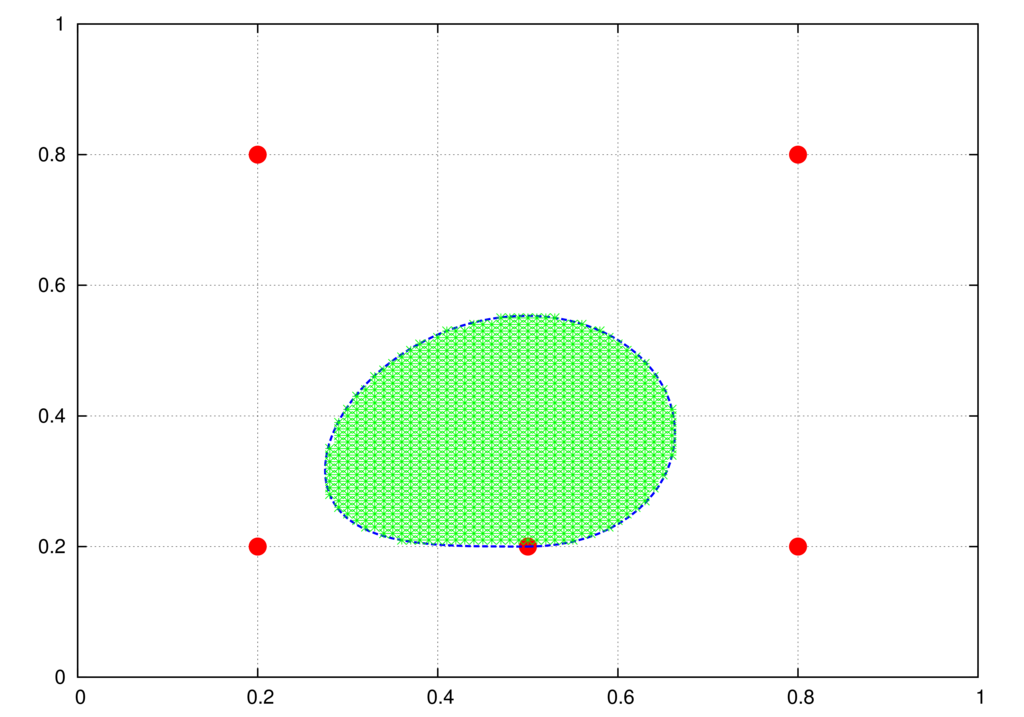}
\includegraphics[width=0.32\textwidth]{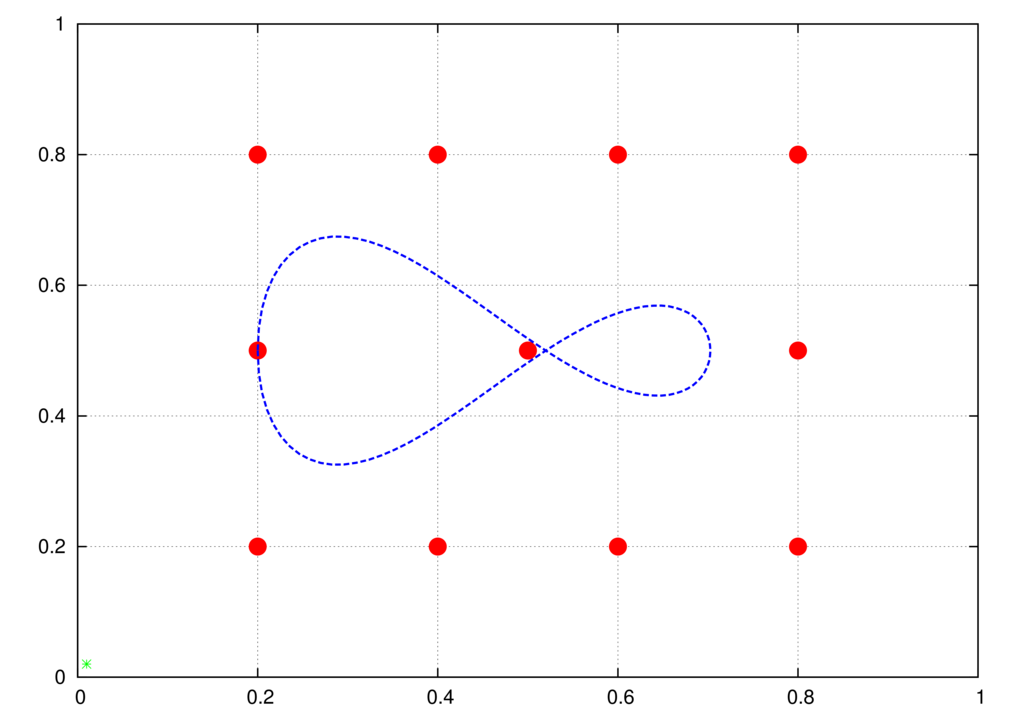}
\includegraphics[width=0.32\textwidth]{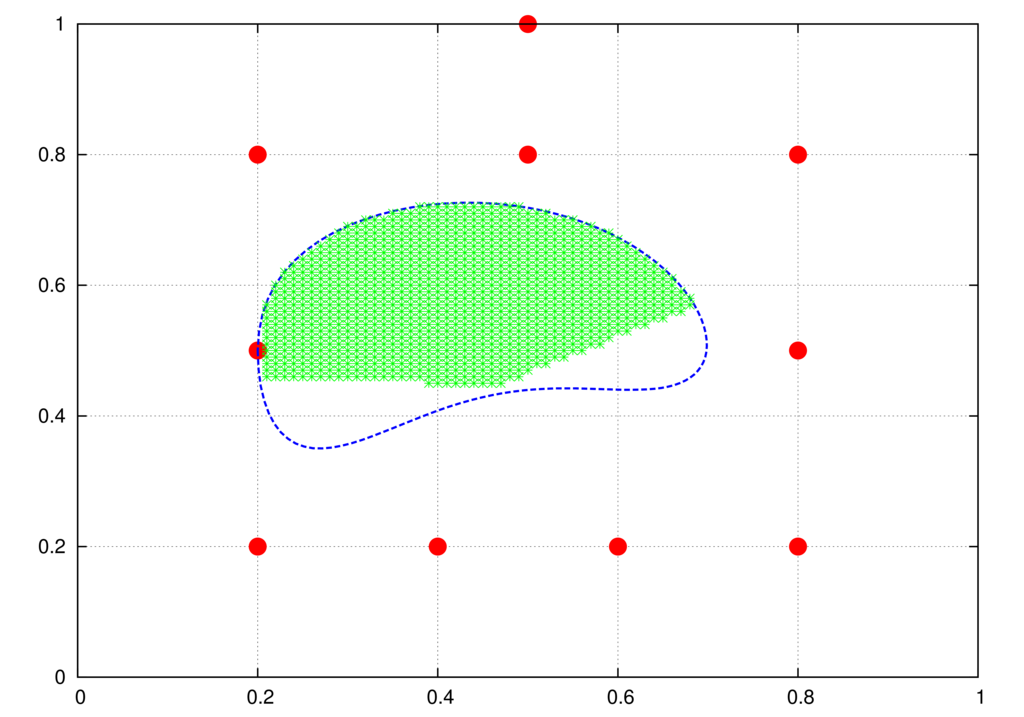}
\caption{Three examples of determination of inner point with respect to a
         closed curve: from left to right, regular curve, self intersecting
         curve, concave curve.}\label{Interni}
\end{center}
\end{figure}

An example is reported in figure \ref{Interni}: three
different situations are depicted (B\'ezier curves of different degree are adopted). 
Red dots are representing the control points of the B\'ezier curve, blue
curve is the resulting B\'ezier curve, green dots are the inner
points as identified by the previously described procedure. On extreme left of figure
\ref{Interni}, the case of a convex shape is presented: all the internal points are
correctly identified. In the middle, the case of a self intersecting curve is represented:
in this case, no internal point is detected, and the obstacle is ignored. On the extreme
right, a concave shape is tested: in this case, not all the inner points are correctly
detected, but the subset of inner points is producing a convex shape. This last quality,
although not respecting the requirement of bijection between shape
and parameters, is appreciated and left unchanged in the parametric model.

In this paper a model with 5 control points will be adopted ($n=4$): being the
contour laying on a 2D space, 6 design variables are required. The objective (cost) function $F:\R^{6}\to\R^+$ to be minimized is finally given by the sequence of the following steps:
\begin{enumerate}
\item Set 3 control points ${\bf P}_{0},{\bf P}_{1},{\bf P}_{2}$;
\item add two further points ${\bf P}_{3}$ ad ${\bf P}_{4}$ to get a closed smooth curve;
\item compute the B\'ezier curve;
\item compute the inner points of the B\'ezier curve;
\item run simulation with the resulting obstacle;
\item get the evacuation time.
\end{enumerate}
\subsection{Optimization algorithm}\label{sec:PSO}

Since we have not a preference for a particular type or class of shapes, we need to investigate
the design space without any prejudice and extensively. As as consequence, a global optimization algorithm is representing the best choice for the problem in hands.

Among the different algorithms already available in literature, the Particle Swarm Optimization
(PSO) represents a really interesting algorithm, whose main qualities are suitable for this
specific application.

First of all, the algorithm is not requiring derivatives of the objective function. The here
adopted mathematical model is not providing derivatives, and its outcome is also quantized
due to the influence of the adopted time step and grid cell dimension: as a consequence, a
continuous variation of a single design variable is not resulting into a continuous variation
of the evacuation time, since the time resolution of the model is connected with the grid cell
dimension. Furthermore, a variation of the obstacle is detected by the mathematical model if
and only if a new grid point is changing his status from internal to external, or vice versa:
as a consequence, the code may result to be dumb with respect to small variations of the
design parameters, if the deformation of the obstacle is not sufficiently large to produce a 
variation of the included or excluded points on the computational grid. Lastly, variation of
the shape of the obstacle may occur in a region not involved by the crowd motion: in this case,
there is no effect on the total evacuation time. For all these reasons, we can expect not to have
a single minimum point, but large plateaux, a finite region sharing the same minimum value, cf.\ \cite[Figs.\ 9 and 14]{cristiani2015SIAP}.

As a consequence of this last feature, the basin of attraction of the global minimum is supposed not to be small. 
PSO is pretty fast in the identification
of the basins of attraction of the objective function, but it is not fast in the identification of
the exact position of the minimum of the selected basin: given the above, this feature is perfectly
compatible with the characteristics of the current objective function.

In PSO, a swarm of elements is distributed onto the design space. An initial speed is assigned
to every element, and the swarm is evolving according to a simple relationship.  
Position of the swarm is updated at each iteration, and the objective function is computed at
the locations identified by the particles. 
Classical formulation of the PSO iteration, originally proposed in \cite{Kennedy}, is provided
in a more general form in \cite{Eberhart}
\begin{equation}
\left\{
\begin{array}{l}
\vv^{k+1}_{i}=\chi \Big( \omega \vv^{k}_{i} +
 c_{1} \left ({\mathbf p}^{k}_{i}-\x^{k}_{i} \right ) +
 c_{2} \left ({\mathbf b}^{k}-\x^{k}_{i} \right ) \Big),
 \\ [2mm]
\x^{k+1}_{i}=\x^{k}_{i}+\vv^{k+1}_{i},
\end{array}
\right.
\end{equation}
where $k$ is the iteration number, ${\bf x}_i$ is the position of the $i^{th}$ particle
of the swarm in the design space, ${\bf v}_i$ is the velocity of the $i^{th}$ particle, ${\bf p}_i$ is the best point ever met by the $i^{th}$ particle in the previous iterations
(corresponding to the personal best value of the objective function), ${\bf b}^k$ is the overall best point ever experienced by the whole swarm, $c_1$ and $c_2$ are two constants
identifying the balance between the local and global phase of the search, $\chi$ is a speed limiter, $\omega$ is the inertia effect.  In the present paper, the following values
are assumed: $c_1=1.494$, $c_2=1.494$, $\chi=1.0$, $\omega=0.729$, while initial position and speed of the initial swarm are assumed according to \cite{peri2015EngComp}.

\REV{From the computational standpoint, the algorithm is intrinsically parallel: in fact, once the
new position of the swarm has been obtained, the computation of the objective function on one
particle is not influencing the computation of the objective function for a different particle.
If we have the availability of one processor for each
particle, the computation of a single iteration is performed in one single shot. In the following, an MPI implementation of the algorithm
will be applied, and the reduction of the computational time is linear with the number of
processors, since the computational effort is mainly related to the numerical model than to the
optimization algorithm.}

In \cite{peri2015EngComp}, a modification of the algorithm for constrained problems is proposed
and tested. The main differences with the original algorithm are:

\begin{itemize}
\item A local search based on a {\it meta-model} \cite{peri2009STR} of the objective function is replacing
      the original PSO step when the overall best particle is not improving further. A
      trust region approach is also implemented in order to guarantee local convergence
      and consistency of the new formulation with the original optimization problem.
\item When two particles are too close to each other, one particle is shifted toward an
      uninvestigated region of the design space. This is here possible because a local
      search phase is integrated into the algorithm, and the local convergence inside a
      basin of attraction is guaranteed. Otherwise, the co-presence of more particles
      in the same basin of attraction is necessary in order to exploit the local search.
\item The search is performed asynchronously until a feasible point is detected. Computation
      of the objective function is performed (possibly in parallel) only once all the
      particles are in a feasible position. Only feasible points are computed, eliminating
      the computation of points that cannot be optimal since constraints are violated. The
      parallel nature of the algorithm is preserved, while useless points are not computed.
\end{itemize}

This new version of the algorithm demonstrates a higher efficiency with respect to the
original formulation when the dimension of the feasible set is reduced by the active constraints. 
The use of a surrogate
model to compute the local behaviour of the objective function may help in reducing the 
small discontinuities induced by the quantized output of the mathematical model.


\section{Preliminary numerical investigation}\label{sec:preliminarynumericalinvestigation}

\subsection{Convergence analysis of the mathematical model}

\begin{figure}[h]
\begin{center}
\begin{overpic}
[width=0.5\textwidth]{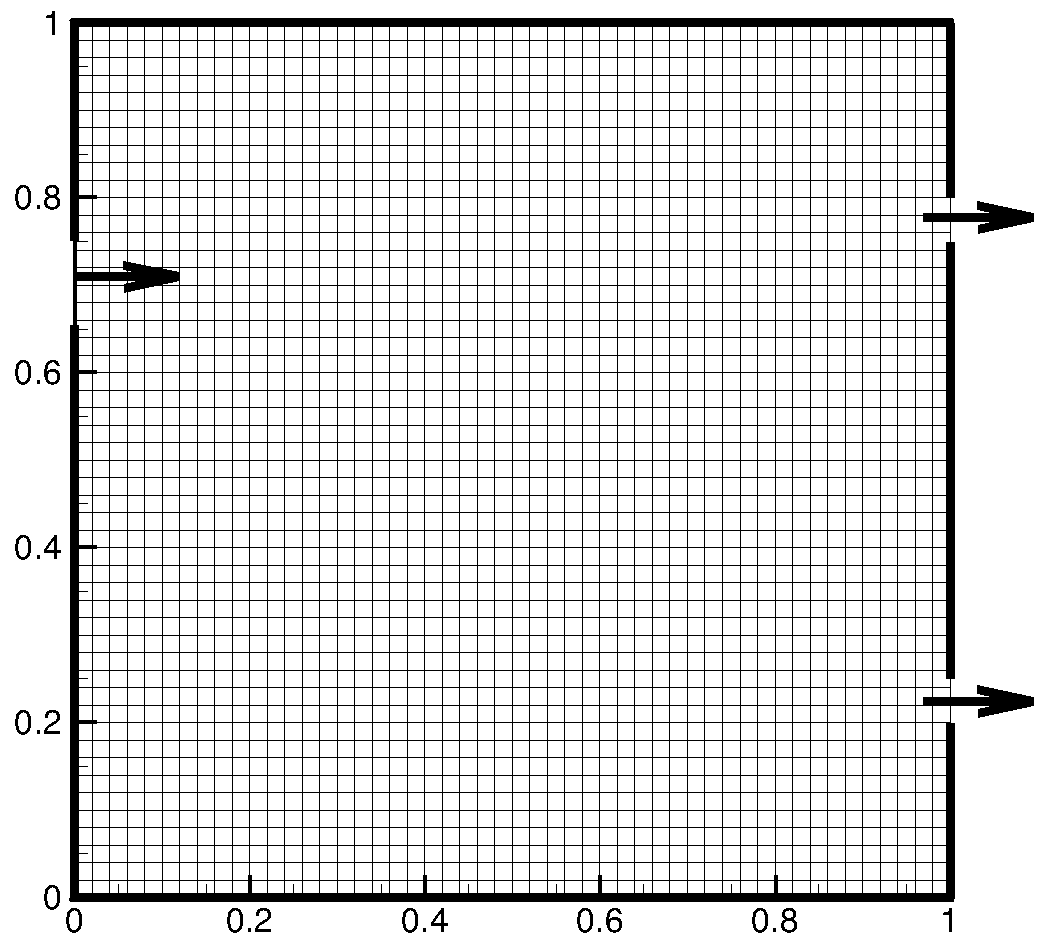}
\put(-2,62){$\mathcal I$}
\put(100,68){$\E$}
\put(100,22){$\E$}
\end{overpic}
\caption{Scheme of the positioning of entrance and exit in the optimization case study. Coarse mesh (50 $\times$ 50) is reported.}
\label{room}
\end{center}
\end{figure}

At the end of the optimization process, the performances of the optimized system
are improved, and the objective function if reduced 
by a quantity which depends on the original configuration, the parameterization scheme,
the number of design variables, etc. However, we have to consider the sensitivity of the
analysis tool we are using to predict the performances and the objective function.
If the improvement is of the same order of magnitude of the accuracy of the adopted
prediction tool, we cannot state that an improvement has been actually achieved, since we are
below the uncertainty level of the model. As a consequence, we need to quantify
the uncertain level of the mathematical model in order to certify the following
optimization results. Following the guidelines in \cite{roache1997}, the Grid-Convergence
Index (GCI) has been computed in order to quantify uncertainties addressed to the
spatial resolution. Three different levels of mesh size have been adopted, doubling
each time the number of grid points on each dimension. The convergence study has
been performed without considering any obstacle in the field. A square room $\Omega$, with one single
entrance $\mathcal I$ and two exits $\E$ has been selected: geometry is shown in figure \ref{room}.
We denote by $\rhoin(t)$ the pedestrian density at entrance. The rest of the boundary $\partial\Omega\backslash(\mathcal I\cup\E)$ is treated as an internal obstacle.
The entrance is nearly in-line with one of the way out, so that the entering crowd is
exploiting a single exit when the room is empty. This situation is going to change
once an obstacle will be added.
Regarding the model parameters, $\Cr$ is set to 11, $r$ is 0.08, $\rhoin$ is 0.75 and pedestrians are entering for a time $0\leq t\leq t_1=1.5$. In this way the room is emptied before the end of the simulation, and the evacuation time can be computed. 
In the following, we will refer to this set-up as ``Test A''.

Results of the convergence analysis are reported in table \ref{tab:GCI}: a non-monotone convergence is observed, with a rate
of convergence of about 3.5 and an uncertainty level largely smaller than 0.5\%. As a consequence,
if the achieved improvement will be higher than 0.5\%, we can certify that an improvement
has been obtained. 
\begin{table}[h!]
\caption{Grid-Convergence study for the adopted mathematical model.
         }\label{tab:GCI}\centering\small
\begin{tabular}{|c|c||l|c|} \hline
Grid level & Evacuation time   & Asymptotic value      & 4.221  \\ \hline
 50 $\times$  50      & 3.942  & Convergence order     & 3.544  \\ \hline
100 $\times$ 100      & 4.257  & Apparent uncertainty   & 0.008  \\ \hline
200 $\times$ 200      & 4.230  & Percentage uncertainty & 0.179  \\ \hline
\end{tabular}
\end{table}

\subsection{Constraints of the optimization process}

When an optimization problem solution is tackled, some further boundary conditions need to be defined,
in order to prevent undesired configurations.
Some geometrical constraints are adopted in order to prevent unusual or unrealistic
situations. First of all, the control points of the B\'ezier curve are limited inside
the volume of the room reduced by three cells along the borders. Since the B\'ezier curve
is included into the convex hull formed by the control polygon, by definition, this guarantees that the obstacle is entirely inside the room. Furthermore, in order not to
create narrow passages between the walls and the obstacle, neither to lock entrance or exits,
a gap has been included.

A second constraint is related to the obstacle's size. It is 
limited in between a minimum and a maximum, measured as a percentage of the room space.
This is because we do not want the obstacle to be too small or too large.

Two further parameters, related to the mathematical model, strongly
influence the final results. The first is the repulsion strength index ($\Cr$),
ruling the interactions among people in the crowd.
The second is the amount of walking people (inflow condition $\rhoin$). 
In fact, it is intuitive to think that, if only few people are passing through the room or they barely affect each other, the obstacle will be disadvantageous. In order to tune these parameters, and also to quantify their effects,
a dedicated sensitivity study has been performed.
\subsection{Sensitivity analysis of the mathematical model}

In order to quantify the effect of the variation of both $\Cr$ and $\rhoin$, a systematic
study has been carried out. 13 different PSO optimization problems have been solved,
covering with regularity the range in between 0.75 and 1.5 for $\rhoin$ and in between 6 and 18 for $\Cr$.
In any case the inflow is active only for a certain time $0\leq t\leq t_1=1.5$, then $\rhoin$ is set to 0 for $t>t_1$. The initial condition is $\bar\rho\equiv0$ (no pedestrians already in the room at the initial time) and $r=0.04$.
In this test, limits for the portion of the room fillable by the obstacle are in between 5 and 25\% of the full room area.

\begin{figure}[h]
\begin{center}
\includegraphics[width=0.9\textwidth]{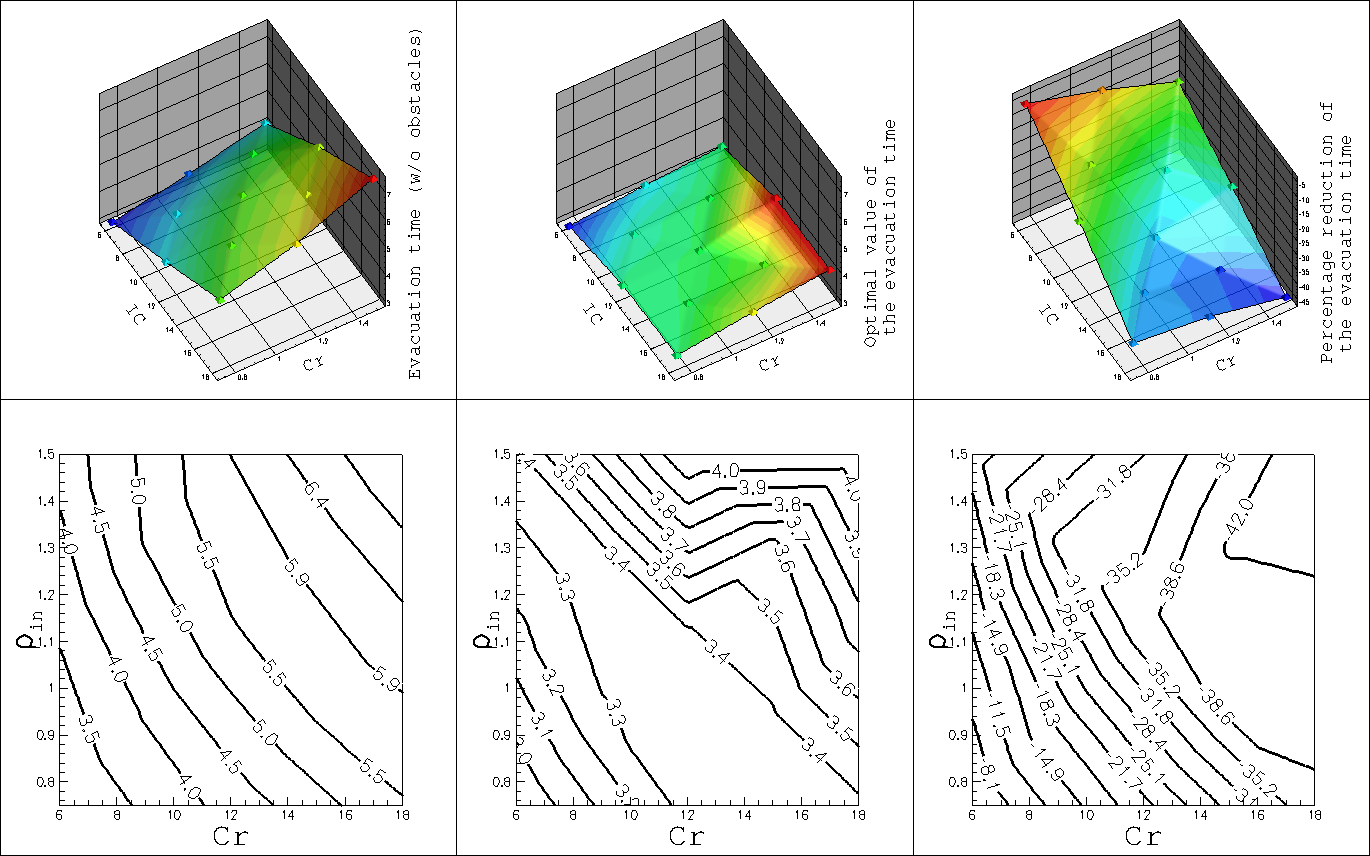}
\caption{Sensitivity analysis of the optimization process with respect to the
         repulsion strength index ($\Cr$) and the inflow density ($\rhoin$).}
\label{SA}
\end{center}
\end{figure}

Since a quite large number of optimization is here performed, the coarse grid (50$\times$50) is adopted.
Results are reported in figure \ref{SA}. In the top part of the picture,
a three-dimensional view and the 13 reference points are reported, while contour levels of the same graph as in the upper part, but with a two-dimensional view, are reported at the bottom part of the picture.
From left to right are reported: the initial value of the evacuation time, without obstacles,
the optimal value of the evacuation time after the insertion of an obstacle, the percentage
improvement. Surfaces are obtained by linear interpolation.

Both $\rhoin$ and $\Cr$ have a strong impact on the final results of the optimization process, with
a higher effect coming from $\Cr$. On top left of figure \ref{SA} the initial evacuation time
is reported: the higher evacuation time is obtained with the higher values of $\rhoin$ and $\Cr$.
The evacuation time is nearly doubled with respect to the condition where both $\rhoin$ and $\Cr$
assume the minimum investigated value. Observing the values on the edges, variation of the
evacuation time with constant $\Cr$ is in between 1.2 and 2.0 time units while $\rhoin$ is passing
from the minimum to the maximum. For the variation of $\Cr$ at constant $\rhoin$ the variation
is in between 1.4 and 3.2 time units. This is demonstrating that, if the adopted limits
for the selected coefficients are representative of the extreme between two limit situations,
the effect of $\Cr$ on the variation of the evacuation time of an empty room is stronger than $\rhoin$. 

On top right, the improvements for the evacuation time obtained adding one obstacles are reported.
Negative values are for improvements, so that the lower the value, the larger the improvement.
Optimization produces a reduction of the evacuation time in between 5 and 45\%, with the
higher values for the highest values of $\rhoin$ and $\Cr$. While for the evacuation time the
trend was not far from being linear, here a more curved profile is observed. 
Anyway, looking at
the extreme values on the edges, improvement pass from 5 and 20\% for the lowest (constant)
value of $\Cr$ when $\rhoin$ is varied, and from 38 to 45\% for the highest value of $\Cr$. This means that
the effects of $\Cr$ and $\rhoin$ are not simply cumulative, but the mixed effect is reducing the strength
of each other. Same situation is observed if $\rhoin$ is fixed and $\Cr$ is varied: for the minimum value
of $\rhoin$, improvement is in between 5 to 38\%, while for the highest $\rhoin$ improvement is 
in between 20 and 45\%. Also in this case, the effect of $\Cr$ is higher than the effect of
$\rhoin$ on the achievable improvements.

\begin{figure}[h]
\begin{center}
\begin{overpic}[width=0.95\textwidth]{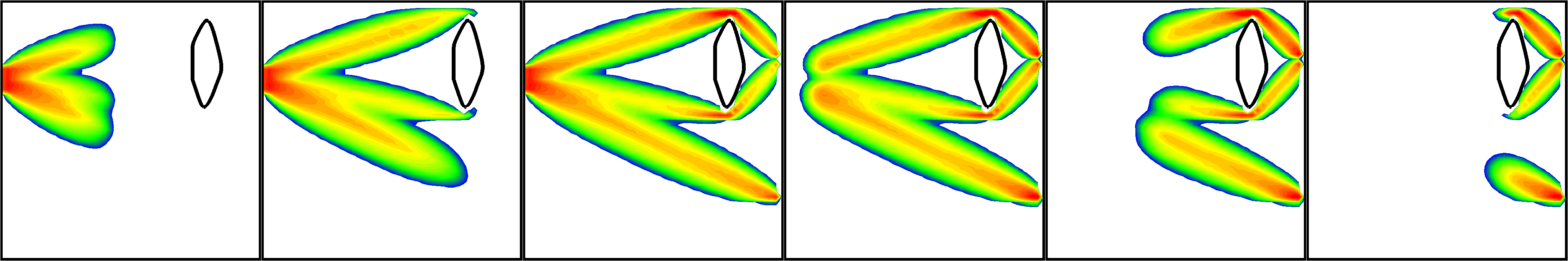}\put(1,1){\tiny $\Cr\!=\!6$} \put(1,3.7){\tiny $\rhoin\!=\!0.75$} \end{overpic}\\ [2mm] 
\begin{overpic}[width=0.95\textwidth]{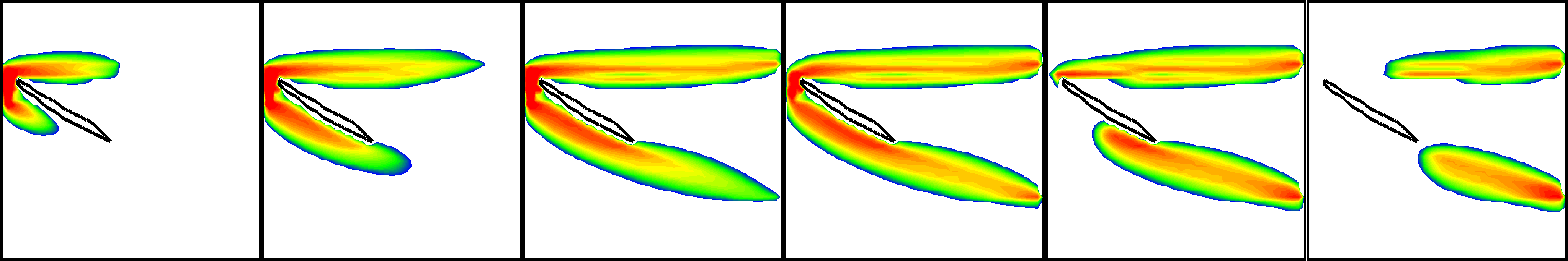}\put(1,1){\tiny $\Cr\!=\!12$}\put(1,3.7){\tiny $\rhoin\!=\!0.75$} \end{overpic}\\ [2mm] 
\begin{overpic}[width=0.95\textwidth]{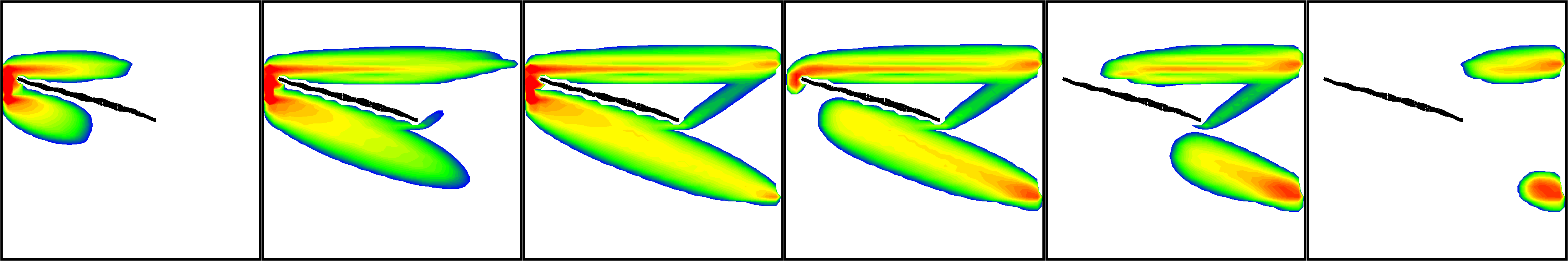}\put(1,1){\tiny $\Cr\!=\!18$}\put(1,3.7){\tiny $\rhoin\!=\!0.75$} \end{overpic}\\ [2mm] 
\begin{overpic}[width=0.95\textwidth]{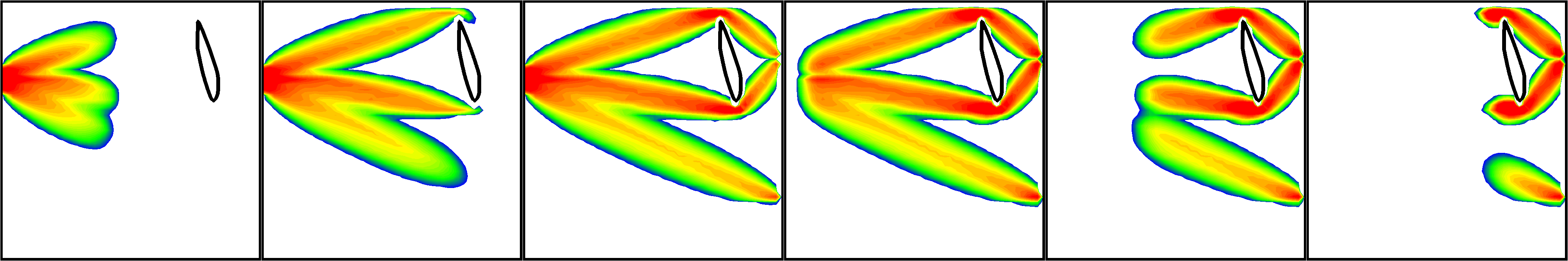}\put(1,1){\tiny $\Cr\!=\!6$} \put(1,3.7){\tiny $\rhoin\!=\!1.125$}\end{overpic}\\ [2mm] 
\begin{overpic}[width=0.95\textwidth]{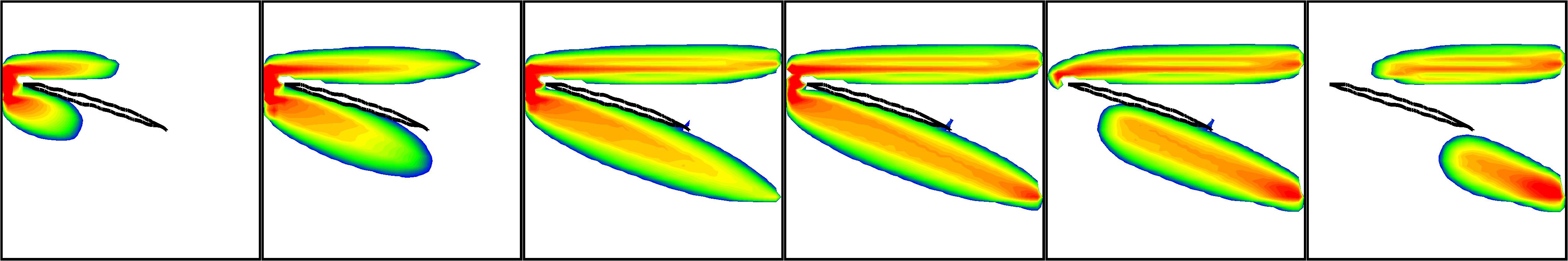}\put(1,1){\tiny $\Cr\!=\!12$}\put(1,3.7){\tiny $\rhoin\!=\!1.125$}\end{overpic}\\ [2mm] 
\begin{overpic}[width=0.95\textwidth]{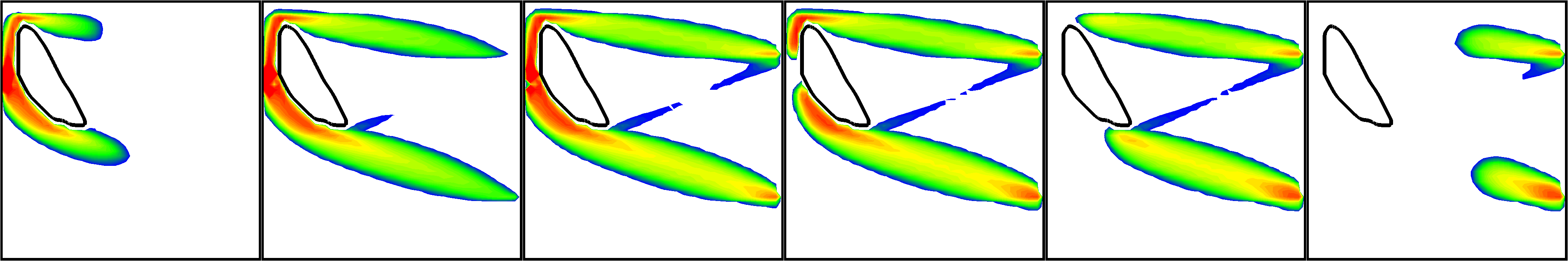}\put(1,1){\tiny $\Cr\!=\!18$}\put(1,3.7){\tiny $\rhoin\!=\!1.125$}\end{overpic}\\ [2mm] 
\begin{overpic}[width=0.95\textwidth]{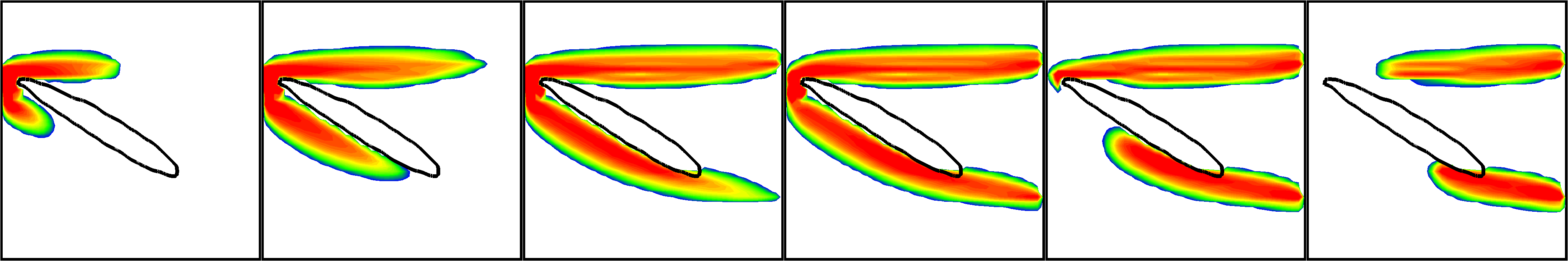}\put(1,1){\tiny $\Cr\!=\!6$} \put(1,3.7){\tiny $\rhoin\!=\!1.5$}  \end{overpic}\\ [2mm] 
\begin{overpic}[width=0.95\textwidth]{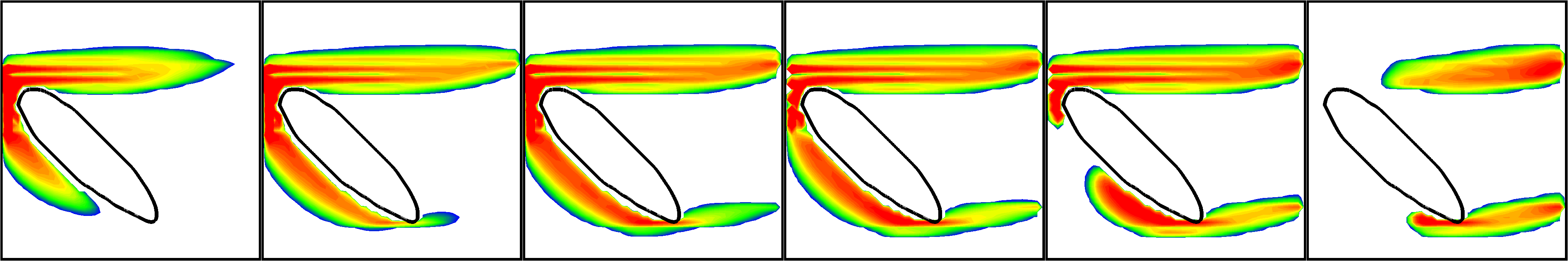}\put(1,1){\tiny $\Cr\!=\!12$}\put(8,8)  {\tiny $\rhoin\!=\!1.5$}  \end{overpic}\\ [2mm] 
\begin{overpic}[width=0.95\textwidth]{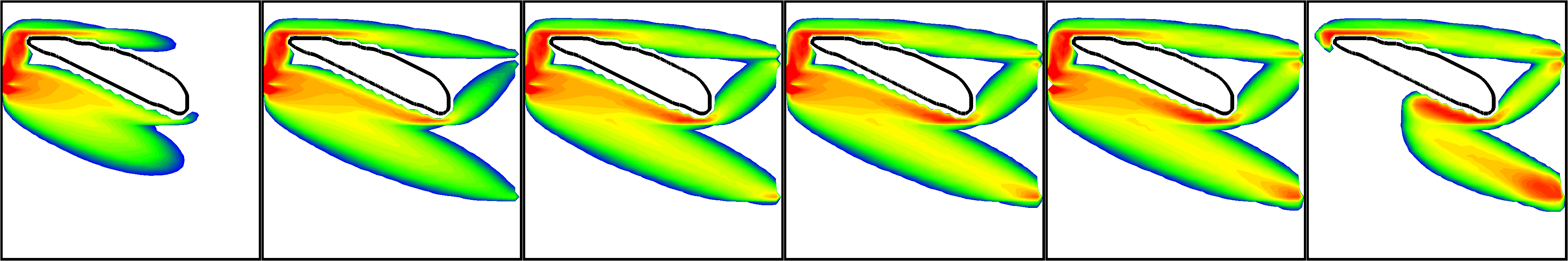}\put(1,1){\tiny $\Cr\!=\!18$}\put(1,3.7){\tiny $\rhoin\!=\!1.5$}  \end{overpic}
\caption{Optimal shape and density evolution for 9 parameter configurations adopted in the
         sensitivity analysis. From left to right, increasing simulation time. 
         See text for information about the visualization of the obstacle.}
\label{SAshape}
\end{center}
\end{figure}

Optimal shapes of the obstacles are reported in figure \ref{SAshape}, together with
a time evolution of the pedestrian density. \REV{Since the convexity of the
obstacle is verified on the computational grid, and a modification of the original shape
is produced once the concave part of the obstacle is removed, the real shape of the obstacle
is somehow influenced by the grid density. In this case, the coarseness of the grid is producing a
stepped contour so, for a better visualization, the curves reported in figure \ref{SAshape} represent the best fit of the external boundary of the obstacle as from the computational grid.}

In figure \ref{SAshape} we can observe how, in all the cases, the crowd
is split by the obstacle so that both the exits are utilized (with no obstacles only the upper exit is used), cf.\ \cite{cristiani2015SIAP}. 
In five cases, crowd is split in three parts: the people running to the bottom exit splits
further in two, and some people is deviating back toward the upper exit. This is creating a
time delay for the people reaching the upper exit, that is beneficial for the overall
evacuation time since prevents congestion and clogging.

\section{Numerical tests}\label{sec:numericaltests}

\subsection{Test A: empty square room with inflow}


Basing on the experience of the previous sensitivity analysis, optimization has
been repeated using the $100\times 100$ grid with the following parameters: $\Cr=11$, $\rhoin=0.75$, $r=0.08$, $t_1=1.5$. 

\subsubsection*{Test A1: influence of the obstacles' shape}
\begin{figure}[h!]
\begin{center}
\includegraphics[width=1.00\textwidth]{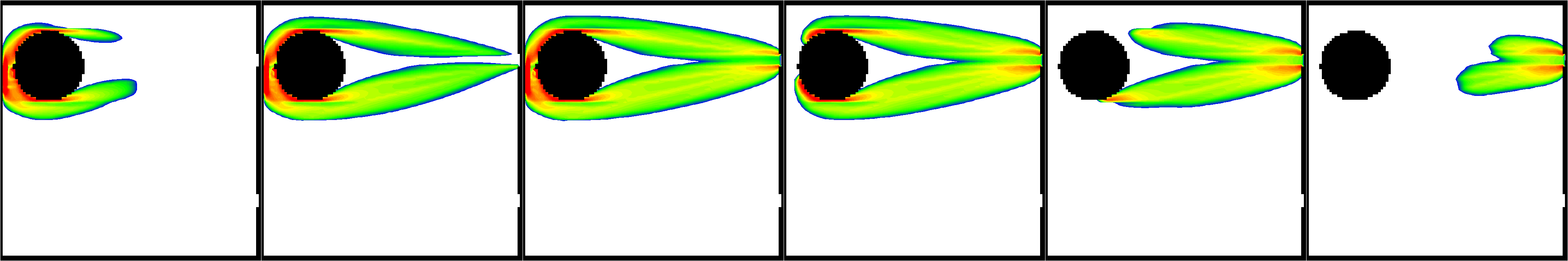}\\ [2mm]
\includegraphics[width=1.00\textwidth]{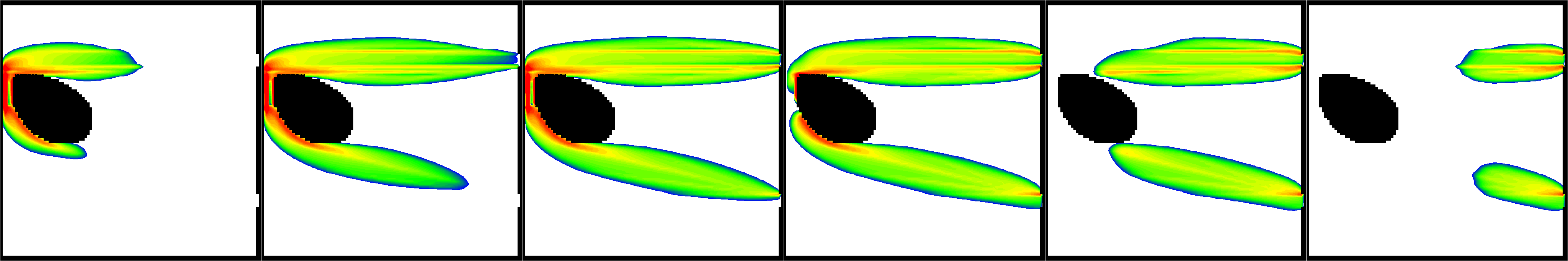}\\ [2mm]
\includegraphics[width=1.00\textwidth]{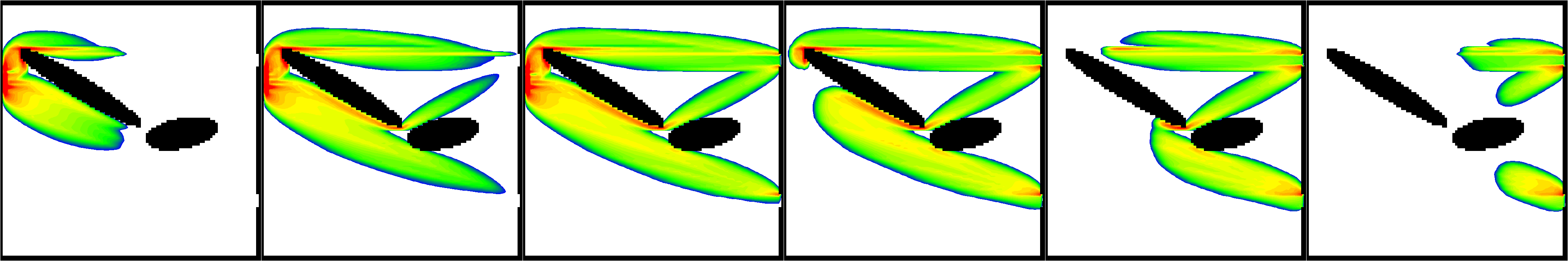}\\ [5mm]
\includegraphics[width=0.9\textwidth]{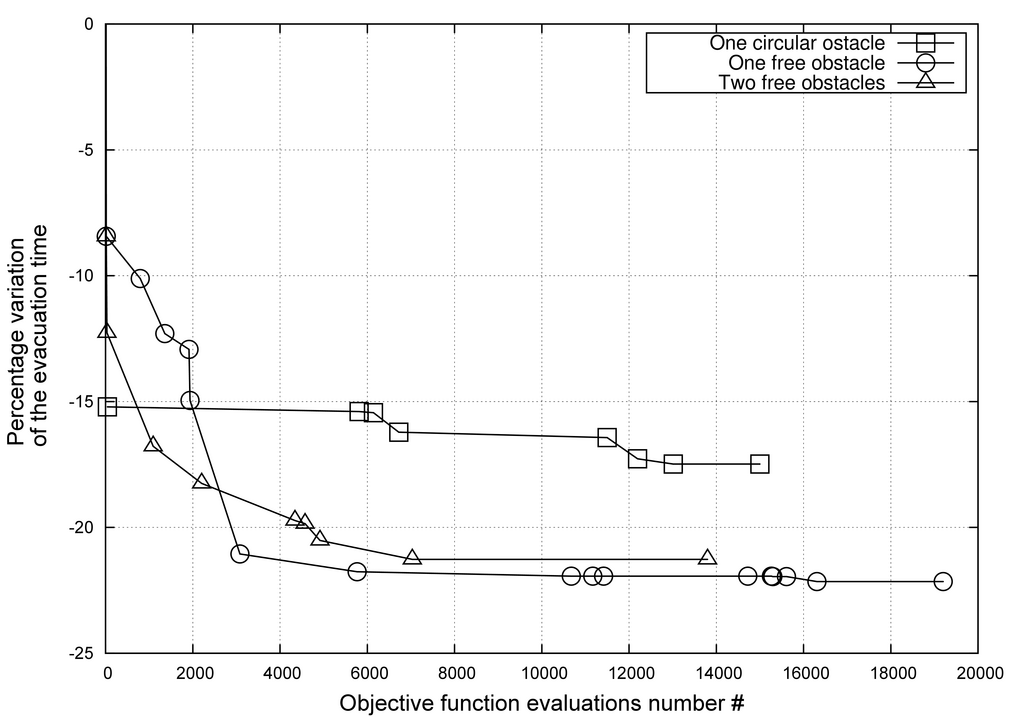}
\caption{Test A1. From top to bottom: evolution of the crowd into the room, observed through 6 different time frames. 
On top, case of one single circular obstacle, in the middle, case of one single free obstacle, on bottom, case of two free obstacles. On extreme bottom, evolution of the best solutions of the previously described problems along the course of the brute force investigation: the three proposed parameterizations are compared.
        }\label{BF}
\end{center}
\end{figure}

A preliminary test has been produced in order to identify
the influence of the shape of the obstacles on the final result. Three different parametrizations are compared:
in the first case, the obstacle is strictly circular, while in the other two examples we have one or two
B\'ezier curves, respectively. Constraint on the dimension of the obstacles are the same: they can occupy a portion
of the full room space in between 5 and 25\%. In order to remove every doubt about the possibility that the optimization
algorithm is not identifying the global minimum, the search of the overall best solution is obtained by a recursive
regular sampling of the design space (brute force investigation). 300 points are dispersed uniformly into the
design space. Once the overall best is identified, a new search is produced, centring 300 new points on the current best point and
spanning a reduced portion of the design space. As each step, the width of the investigated area is reduced by $1/13$ while the number of the sampling points is preserved: 50 iterations have been attempted for the single circular
obstacle (15000 objective function evaluations), 64 for the single B\'ezier obstacle (19200 objective function evaluations) and 46 for the double B\'ezier obstacles (13800 objective function evaluations). The results are reported in
figure \ref{BF}.

From figure \ref{BF}, we can observe how the possibility to have a free shape for the obstacle is very important. 
In fact, the configuration with a single B\'ezier curve outperforms the single circle of about 5 percentage points. 
On the other side, there is not
a substantial difference between a single obstacle and a couple of obstacles, since the final result produced by the two
configurations is comparable. It must to be stressed that the number of design variables is doubled passing from one to two
obstacles, so that the dimension of the investigated space, and the complexity of the problem, are largely increased.

\medskip

It must be stressed that the case with one obstacle is not completely contained into the case with two obstacles:
in fact, due to the presence of a box constraint on the space occupied by the obstacles, the maximum dimension of the bigger
obstacle is reduced when more than a single obstacle is modelled, because a portion of the space is already occupied by other
obstacles. Looking from the perspective of the single obstacle, the constraint limits are changed as a consequence of the
presence of the other obstacles: each obstacle can be smaller than the lower limit, since the other obstacles will contribute
to the final value of the occupied space, but it will be also smaller than the prescribed upper limit, in order to leave some space
to the other obstacles. This situation is reflected in the final results reported in figure \ref{BF}: the increase of the number
of the obstacles is not causing a further reduction of the evacuation time. Probably a modification on the applied constraints
would be beneficial in this sense.

\subsubsection*{Test A2: influence of the obstacles' size}
Once the importance of the shape of the obstacle has been evidenced, two different optimizations, using two different sets of
constraints, have been produced. One or two obstacles are considered, and the occupation of the room by the obstacles is
constrained in between 10 and 90\%, or in between 10 and 20\%. As a consequence, a total of 4 different optimization problems
have been solved.

\begin{figure}[h!]
\begin{center}
\includegraphics[width=1.00\textwidth]{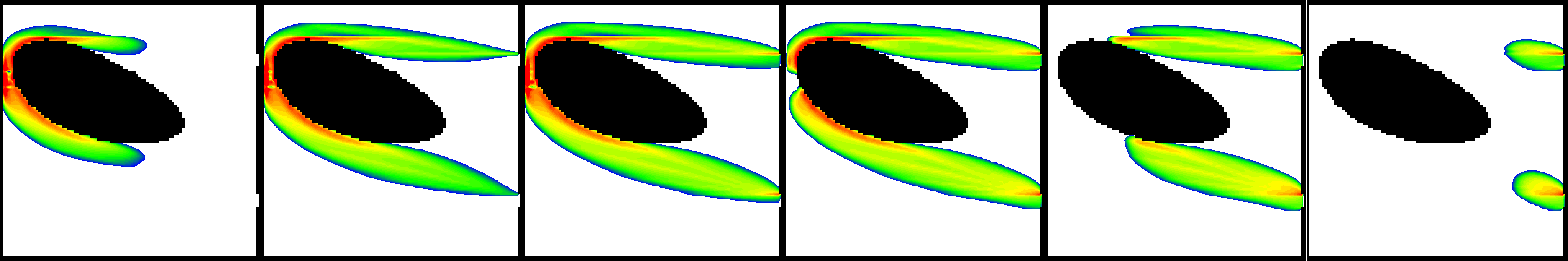} \\ [2mm]
\includegraphics[width=1.00\textwidth]{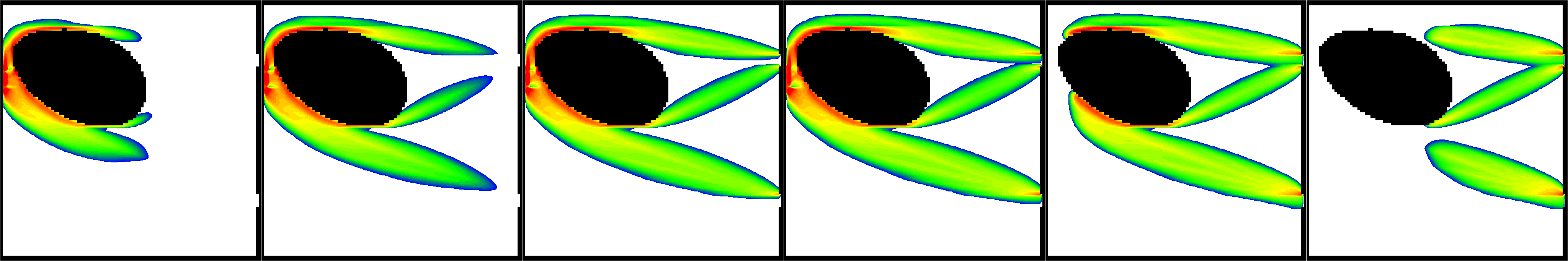} \\ [2mm]
\includegraphics[width=1.00\textwidth]{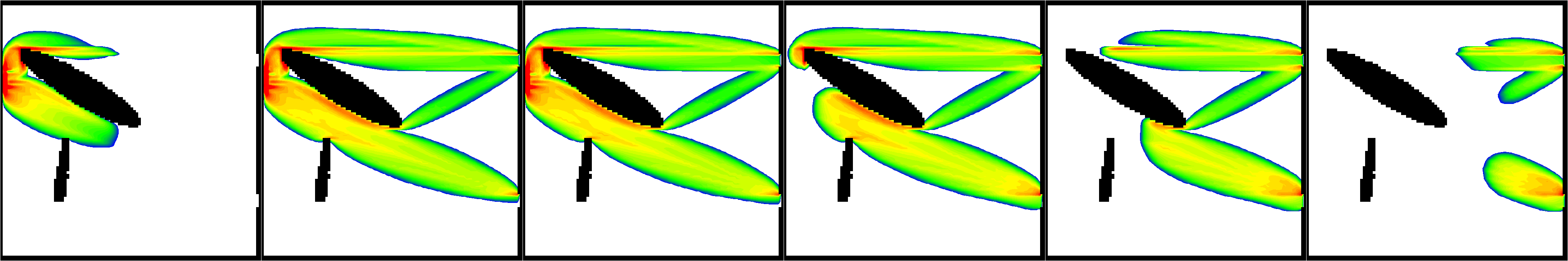} \\ [2mm]
\includegraphics[width=1.00\textwidth]{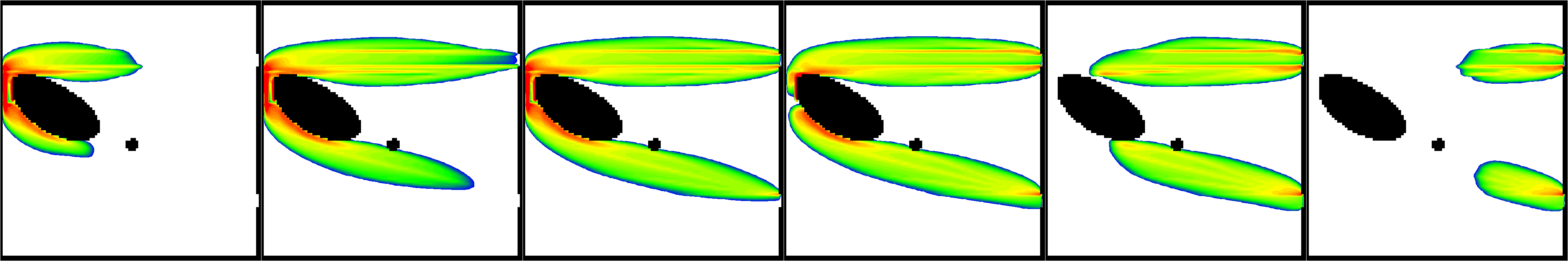} \\ [2mm]
\includegraphics[width=0.9\textwidth]{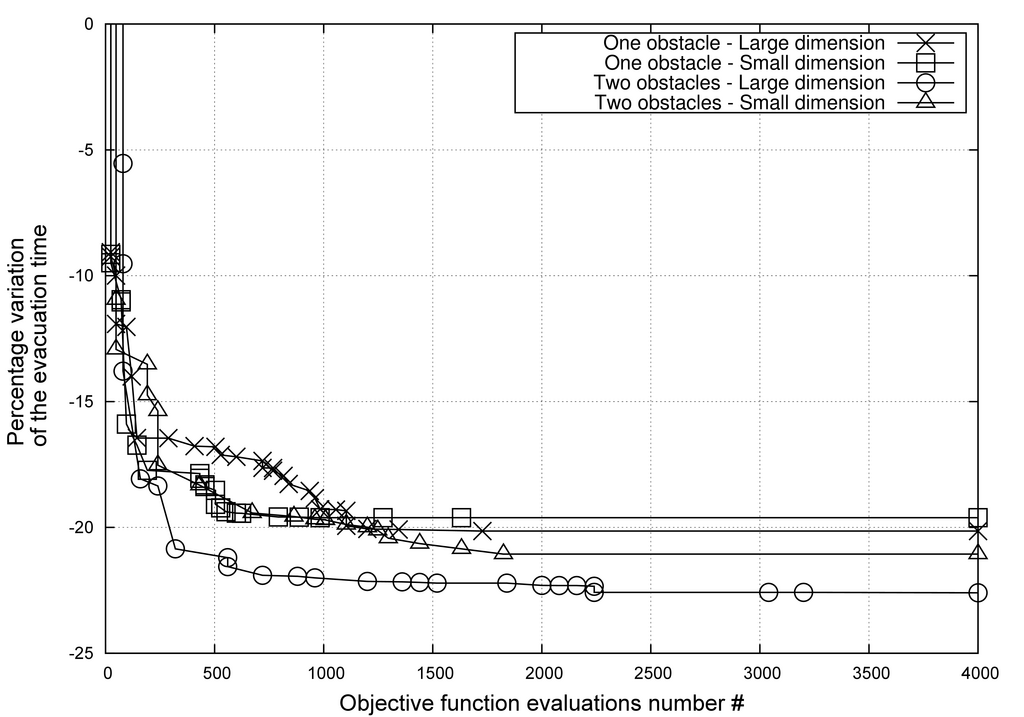}
\caption{Test A2. From top to bottom: evolution of the crowd into the room, observed through 6 different time frames.
         On top, case of one single free obstacle with tight
         constraints on the dimension of the obstacle with respect to the dimension of the room. Second picture
         is representing the same case, but with larger obstacles allowed. Third and fourth picture replicate the test with
         two obstacles. 
         On extreme bottom, the percentage reduction of the evacuation time with respect to the empty room as a function of the iterations of the optimization algorithm.
        }\label{StanzaOpt_e_FlowOptStanza}
\end{center}
\end{figure}

In figure \ref{StanzaOpt_e_FlowOptStanza}, the results produced by the PSO algorithm are presented. Convergence of the algorithm to the optimal solution is evidenced. 
The increase of the number of obstacles is producing, in this case, larger improvement on the final evacuation time: this is probably
due to the relaxation on the size constraint limitation.
The obstacles are dividing the crowd in two or three groups, inducing a time delay in the arrival of the groups to
the exit, with a better exploiting of the available exits. 
The effect of the obstacles, producing a complete use of the available exits and a time shift for the different parts of the crowd, is well documented. 
An example of time shift induction is represented by the fourth case in figure \ref{StanzaOpt_e_FlowOptStanza}: here a very small obstacle is placed in
the lower part of the room. The crowd is marginally touching the obstacle, but this is producing a local density increase in the vicinity
of the small obstacle, so that some of the pedestrians travelling that direction are slowed down, with a subsequent reduction of the
exit crowding and an increase of the flow speed at the exit.

The use of multiple obstacles is obviously increasing the number of parameters, increasing accordingly the dimension of the investigated space:
this is a great difficulty for the optimization algorithm, because the associated problem have a larger complexity. Furthermore, the use of a second
obstacle put some further constraints on the optimization problem to be solved, since the second obstacle, in order to be useful, must lay outside the room
space already occupied by the first obstacle. This is excluding some portions of the parameter space from the feasible space, so that the solution
of the optimization problem is even more difficult to solve. The use of the inner point variant of the PSO algorithm \cite{peri2015EngComp} is alleviating this problem,
and the final solution is absolutely similar to the solution obtained by using a brute force approach, using less than $1/5$ of the objective function evaluations.
\subsection{Test B: partly filled room with no inflow}

As a second test case, we consider a square room with a single exit on the right wall and no entrance. At the initial time, the room is partly filled by steady people in the left half part ($\bar\rho=0.75$ in $[0.1,0.5]\times[0.2,0.8]$). 
In this case, the flow of the people toward the exit is perfectly symmetric, and the optimal positioning 
of the obstacles is expected to break this symmetry in order to facilitate the room evacuation. The other parameters are the same as Test A: $\Cr=11$, $r=0.08$.

Due to the extreme simplicity of the geometry, the use of a single obstacle will probably not be sufficient to produce
substantial improvements. 
On the other hand, the raise of the number of obstacles implies a huge increase of the complexity of the optimization problem, firstly as a consequence of the
higher space dimension, but also due to the introduction of some implicit supplementary constraints, that could have a great
impact on the shape of the feasible set associated with the resulting optimization problem. In fact, the use of multiple free
obstacles induces the necessity to discharge all the solution for which the two obstacles are overlapping: otherwise, this will
cause large zones of the feasible space with constant value of the objective function, because one big obstacle can hide smaller
obstacles without causing a change in the configuration of the room to be considered. In other terms, once one big obstacle is defined,
fixing the related parameters, all the other parameters are irrelevant as soon as they are not able to draw an obstacle laying outside
the first obstacle: consequently, the value of the objective function is piecewise constant.

\subsubsection*{Test B1: circular obstacles}
As a preliminary check, in order also to gain some insights into the problem, a simplified parameterization is adopted, in order
to identify some promising configurations to be further investigated. Here a fixed number of circular obstacles are disposed into
the room. Diameter of the circles is fixed, equal to $\frac{0.5}{N}-0.01$ (if $N$ is the total number of circles - diameter is
decreased by 0.01 in order not to have overlaps between the circles), as well as the longitudinal position of the center of each
circle (defined as $0.5 + \frac{0.5}{N} + \frac{i-1}{N}$, where $i$ is indicating the ordinal circle number, from left to right),
so that, if the circles are aligned, no overlap is obtained while the circles are covering nearly all the empty space. As a
consequence, we have here only $N$ design variables, with a feasible set extremely simple, since each variable is limited in a
prescribed closed interval, without any further limitations.

\begin{figure}[h]
\begin{center}
\includegraphics[width=1.00\textwidth]{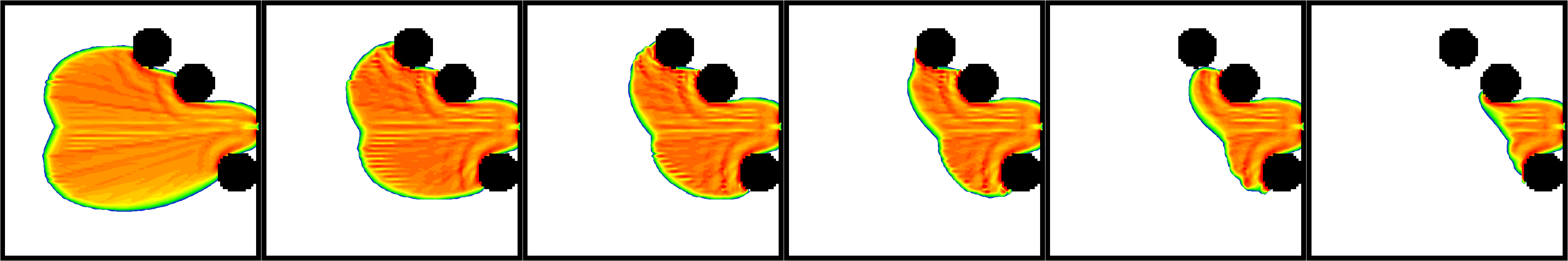}\\ [2mm]
\includegraphics[width=1.00\textwidth]{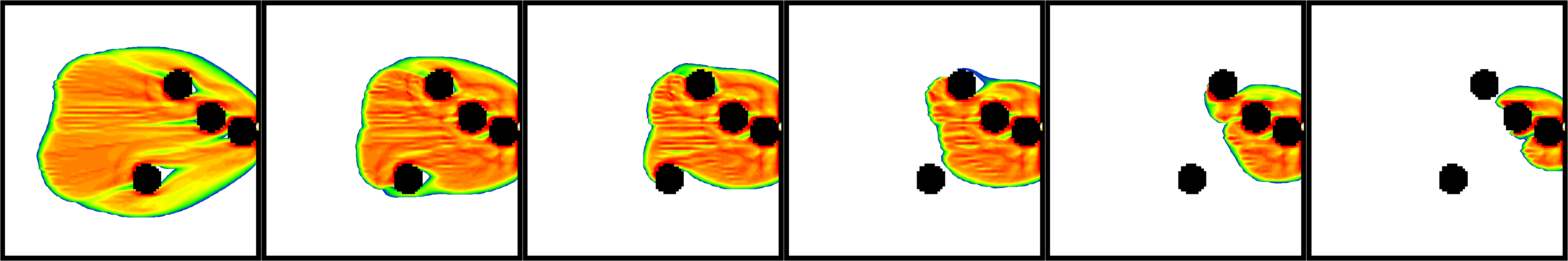}\\ [2mm]
\includegraphics[width=1.00\textwidth]{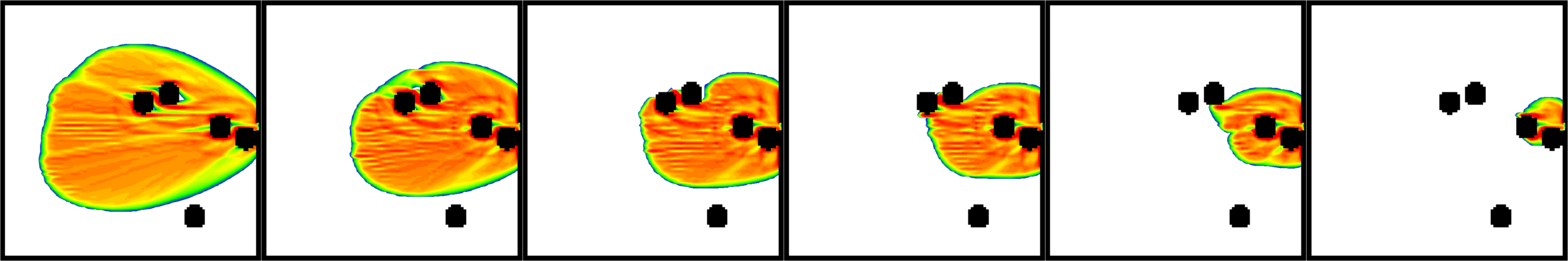}\\ [2mm]
\includegraphics[width=1.00\textwidth]{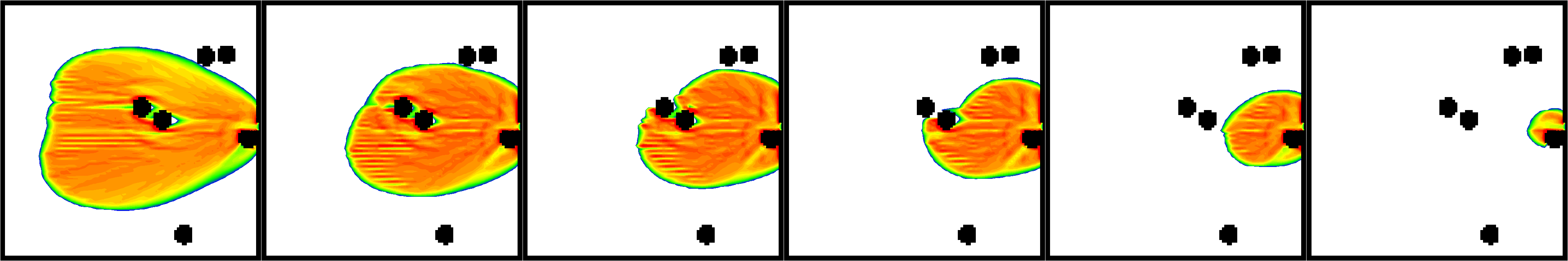}\\ [2mm]
\includegraphics[width=1.00\textwidth]{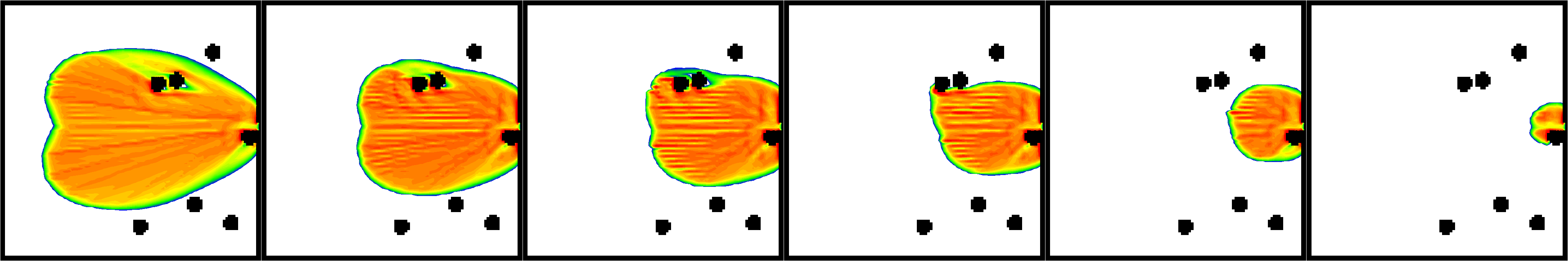}\\ [2mm]
\includegraphics[width=0.84\textwidth]{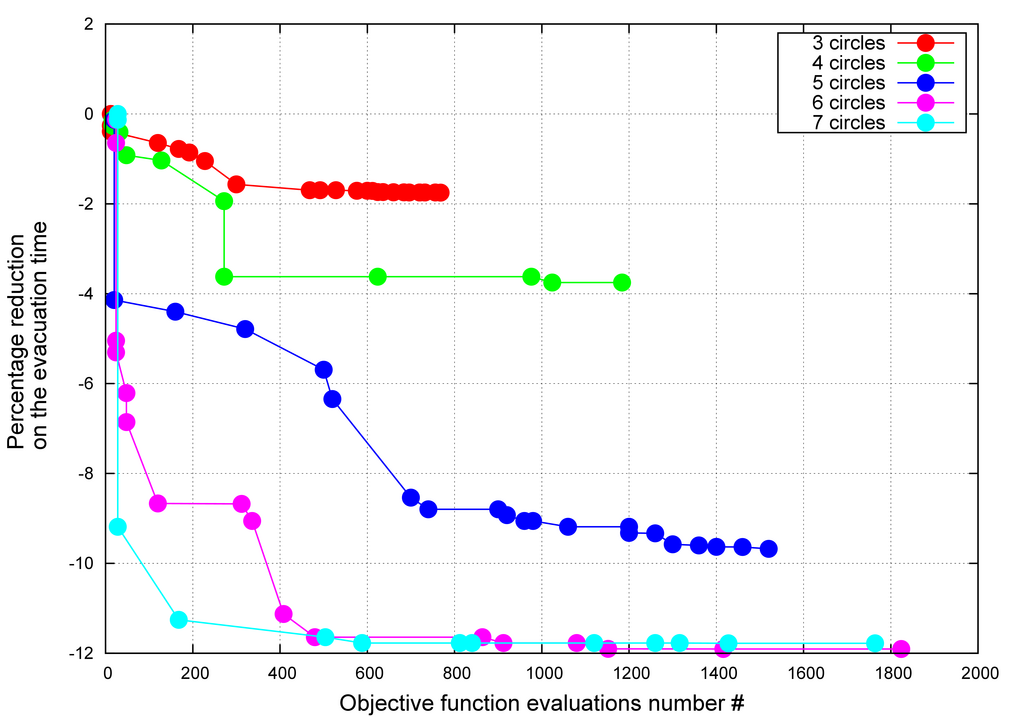}\\ [2mm]
\caption{Test B1. From top to bottom: optimal solution of the half-empty-room problem by using an increasing number of
         staggered circular obstacles. On extreme bottom, the percentage improvements obtained by using the five different
         formulations.
        }\label{HalfRoomStaggered}
\end{center}
\end{figure}

Five different optimization problems have been solved, ranging from three to seven circles. Results are presented in figure
\ref{HalfRoomStaggered}. A 2\% decrease of the evacuation time is obtained with three circles, and this gain is increasing up
to the case with six circles, where the decrease is of about 12\%. No further improvements are observed by passing from six
to seven circles.

Looking at the optimal solutions, and in particular the path of the pedestrians moving from the half left part of the room to
the exit, we can see that while in the case of three and four circles all the obstacles play an active role in the deviation
of the pedestrian flow, in the other cases some of the obstacles are not active, laying in a region never touched by the
passage of the pedestrians. In particular, in the case of six circles, only three obstacles are active. Furthermore, we can see
how all the optimal solutions have one obstacle in vicinity of the exit. This results are suggesting some main conclusions (cf.\ \cite[\S 6.3]{albi1504.04064}):
\begin{itemize}
\item One obstacle may stay close to the exit;
\item small obstacles are needed;
\item three obstacles may be sufficient to guarantee a substantial reduction of the evacuation time.
\end{itemize}

\subsubsection*{Test B2: free-shaped obstacles}
\begin{figure}[h]
\begin{center}
\includegraphics[width=1.00\textwidth]{figs/best-stag-6.png}\\ [2mm]
\includegraphics[width=1.00\textwidth]{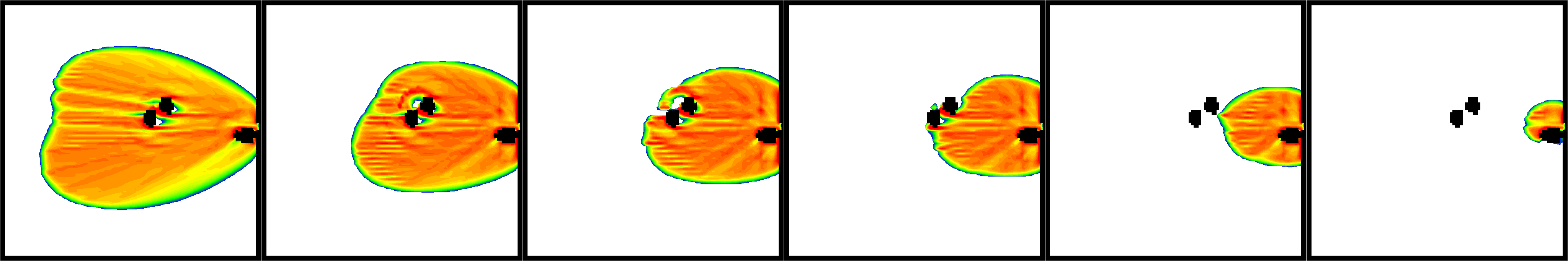}\\ [2mm]
\includegraphics[width=0.9\textwidth]{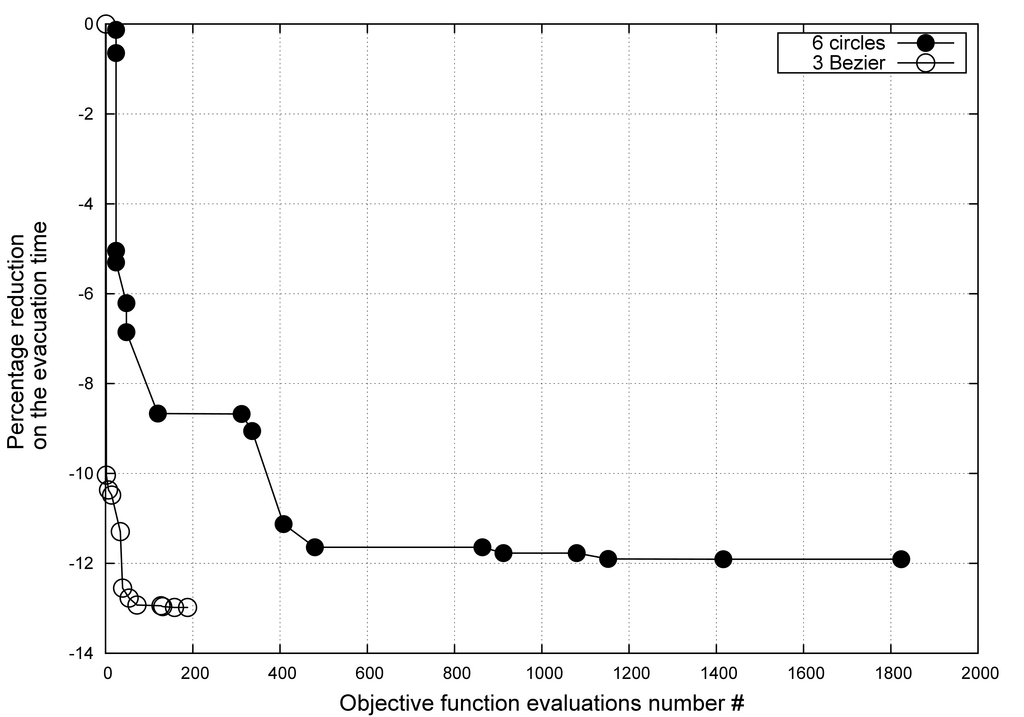}
\caption{Test B2. From top to bottom: optimal solution of the half-empty-room problem by using six staggered circular obstacles
         and three obstacles with unprescribed shape. On extreme bottom, the percentage improvements obtained by using the two different
         approaches.
        }\label{HalfRoomBezier}
\end{center}
\end{figure}

Starting from a configuration similar to the optimal one in the case of six obstacles, a local optimization has been produced,
utilizing three free obstacles only, obtained by using three different B\'ezier curves, as in the case of the empty room: 18
parameters are needed. A pattern search algorithm has been applied here \cite{kolda2004}: this is because we are looking for
a solution in the vicinity of the initial solution, and a local optimization algorithm is able to perform this search much more
efficiently than a global optimization algorithm. Results are presented in figure \ref{HalfRoomBezier}. Here the total number
of solutions required in order to gain convergence is evident: pattern search is requiring one tenth of the iterations, but
it is obviously detecting a local minimum.

A further reduction of about 1 percentage point is obtained, and the final pattern of the obstacles is not far from the starting one.
The preservation of the dimensions of the obstacles is substantially confirming that, in the previous test, the dimension of the
obstacles represents probably the most important parameter, and the increase from six to seven obstacles was reducing too much
the dimension of the obstacles, so that no further improvement was obtained. On the other hand, the possibility of a free-shaped
obstacle is providing a further reduction of the evacuation time, confirming one of the outcome of the previous investigations, that
is, the shape of the obstacle represents an important element for the minimization of the evacuation time.

\section*{\REV{Conclusions, insights, and future work}}

The usefulness of the interposition of one or more free-shaped obstacles between a group of evacuating pedestrians and the exits of a room has been investigated in this paper. Despite the difficulty of the optimization problem, mainly related to the huge number of the possible shapes in comparison with the small number of feasible solutions, the obtained reduction of the total evacuation time is remarkable. Furthermore, the use of an optimization algorithm provides a significant reduction of the total number of configurations to be analyzed before the identification of the optimal solution, increasing the efficiency of the search phase:
a reduction greater than the 80\% is obtained if we compare a systematic/extensive investigation of the design space with the use of PSO. This study is putting on the spotlight the role of the  optimization techniques in the reduction of the evacuation time. 

The importance of local peculiarities in the shape of the obstacles has been also highlighted by comparing free shape with circular shape: the final gain is also influenced by small details, and the advantage may be reduced by an incorrect shaping. This kind of result can be hardly obtained by classical design activity.

\medskip

The reliability of the obtained results is questionable unless a complete validation of the mathematical model against experimental data is performed. It is important to recall that our modelling assumptions could be far from being true, since we are considering here pedestrians who are familiar with the environment, follow the shortest path and do not forecast the behaviour of the group mates (see Remark \ref{rem:limitations}). Moreover, we are not taking into account important crowd features (e.g.\ age, emotional state, family groups, social influence, etc.), and we are using a first-order model which is not able, by definition, to describe acceleration and pressure effects. 
Nevertheless, the reduction of the evacuation time achieved in our simulations is so significant that can be hardly nullified modifying some surrounding circumstances. In particular, one of the effects of well-placed and well-shaped obstacles is to split the crowd in order to make all the exists equally used and/or alleviate congestion by delaying a part of the crowd; moreover, once obstacles are placed near the exits, they can also break the symmetries and avoid deadlocks. We think that these features are sufficiently independent of simulation parameters and crowd behaviour. Let us also note that, in the case of unknown environments, obstacles can address people towards unexplored areas, thus making visible some exits.

Finally, let us note that the 
application of the same techniques to more realistic situations, where the complexity of the shape of the evacuating room is possibly higher, may result in a further confirmation of the good results presented in this paper.
This mostly because the optimization algorithm may help in detecting the most favourable solution in a framework in which the high degree of complexity is hiding the most convenient configuration.

The study presented in this paper can help engineers and architects to design new spaces or increase safety of the existing ones. Obviously, the impact of the new architectural elements must verified by suitable field-test in normal and emergency situations, and a nontrivial problem can arise if a given design turns out to be advantageous in some situation and disadvantageous in another one.


\medskip

The sensitivity of some parameters on the final evacuation time, observed in the paper, suggests
to take into account their variability. For example, the amount/position of people occupying the
room is subject to variability in practice, as well as the incoming crowd density. The intrinsic
stochasticity of these variables may represents an interesting parameter for a different case study.

\subsection*{Acknowledgements}
The authors thank Andrea Tosin for helping to create the model described in section \ref{sec:model}.
\bibliographystyle{plain}
\bibliography{biblio}

\newpage

\end{document}